\documentclass[12pt,twoside]{article}
\usepackage{times}
\usepackage{mathrsfs}
\usepackage{amsmath,amssymb}
\usepackage{color}
\usepackage{fancyhdr}
\allowdisplaybreaks

 \def \nn{\nonumber}

\newcommand{\pf}{\noindent {\bf Proof. \hspace{2mm}}}
\newcommand{\ef}{ \hfill $ \Box $ \vskip 3mm}
\newcommand{\be}{\begin{equation}}
\newcommand{\ee}{\end{equation}}
\newcommand{\bea}{\begin{eqnarray}}
\newcommand{\eea}{\end{eqnarray}}

\newcommand{\bR}{{\mathbb R}}
\newcommand{\bN}{{\mathbb N}}
\newcommand{\bZ}{{\mathbb Z}}

\def\R {\Bbb R}

\def\p{\partial}
\def\la{\lambda}
\def\La{\Lambda}
\def\al{\alpha}
\def\t{\tilde}
\def\q{\quad}
\def\th{\theta}

\def\dl{\delta}

\def\lt{\left}
\def\rt{\right}

\def\dl{\delta}
\def\Dl{\Delta}
\def\i{\infty}

\def \ls{\lesssim}
\def\p{\partial}
\def\f{\frac}
\def\na{\nabla}
\def\al{\alpha}

\def\t{\tilde}

\def\O{\Omega}

\def\q{\quad}
\def\qq{\qquad}

\def\B{\dot{\mathcal{B}}}
\def\E{\mathcal{E}}

\def\dB{\dot{B}}

\def\dD{\dot{\Delta}}

\def\mE{\mathcal{E}}
 \allowdisplaybreaks
\def\mG{\mathcal{G}}
\def\mP{\mathcal{P}}
\def\mvP{\mathcal{P}^{\perp}}
\renewcommand{\div}{\mbox{\rm div}\;\!}
\def\Id{{\rm Id}\,}

\begin{document}
\let\oldsection\section
\renewcommand\section{\setcounter{equation}{0}\oldsection}
\renewcommand\thesection{\arabic{section}}
\renewcommand\theequation{\thesection.\arabic{equation}}
\newtheorem{claim}{\noindent Claim}[section]
\newtheorem{theorem}{\noindent Theorem}[section]
\newtheorem{lemma}{\noindent Lemma}[section]
\newtheorem{proposition}{\noindent Proposition}[section]
\newtheorem{definition}{\noindent Definition}[section]
\newtheorem{remark}{\noindent Remark}[section]
\newtheorem{corollary}{\noindent Corollary}[section]
\newtheorem{example}{\noindent Example}[section]

\title{Global existence in critical spaces for non Newtonian compressible viscoelastic flows}

\author{Xinghong Pan,\ Jiang Xu \& Yi Zhu
\vspace{0.5cm}
}

\date{}

\maketitle

\centerline {\bf Abstract} \vskip 0.3 true cm
We are interested in the multi-dimentional compressible viscoelastic flows of Oldroyd type, which is one of non-Newtonian fluids exhibiting the elastic behavior. In order to capture the damping effect of the additional deformation tensor, to the best of our knowledge, the ``div-curl" structural condition plays a key role in previous efforts. Our aim of this paper is to remove the structural condition and prove a global existence of strong solutions to compressible viscoelastic flows in critical spaces. The new ingredient lies in the introduction of effective flux $(\theta,\mathcal{G})$, which enables us to capture the dissipation arising from \textit{combination} of density and deformation tensor. In absence of compatible conditions, the partial dissipation is found in non-Newtonian compressible fluids, which is weaker than that of usual Navier-Stokes equations.

\vskip 0.3 true cm

\vskip 0.3 true cm


\section{Introduction}
In the Eulerian description,  a general compressible fluid evolving in some open set $\Omega$ of  $\R^n$ is characterized
at every  material point~$x$ in~$\Omega$ and time  $t\in\R$
by its  {\it velocity field} $u= u(t,x)\in\R^n,$   {\it  density} $\varrho=\varrho(t,x)\in\R_+,$
{\it pressure} $\Pi=\Pi(t,x)\in\R$. In the absence of external forces and heat diffusion, those physical quantities are governed by
\begin{itemize}
\item The mass conservation:
$$
\partial_t\rho+\div(\rho u)=0.
$$
\item The momentum conservation:
$$
\partial_t(\rho u)+\div(\rho u\otimes u)=\div\mathbb{S}-\nabla\Pi.
$$
\end{itemize}
In the regime of \emph{Newtonian fluids}, $\mathbb{S}$ stands for the \emph{viscous stress tensor}, which is given by
$$
\mathbb{S}\triangleq \lambda\div u\,\Id+2\mu D(u).
$$
Here $\lambda$ and $\mu$ are the {\it viscosity coefficients}
and $D( u)\triangleq\frac12(\nabla u+{}^T\!\nabla u)$ is the {\it deformation tensor.} So one can get the baratropic Navier-Stokes equations of compressible fluids:
\begin{equation}\label{R-E1}
\left\{\begin{array}{l}
\partial_t\rho+\div(\rho u)=0,\\[1ex]
\partial_t(\rho u)+\div(\rho u\otimes u)-\div\bigl(2\mu D(u)+\lambda\,\div u\, \Id\bigr)+\nabla\Pi=0.
\end{array}\right.
\end{equation}

In the last several decades, there have been many attempts to capture different phenomena for \emph{non-Newtonian fluids} such as those in the Ericksen-Rivlin models, the high-grade fluid models, the Ladyzhenskaya models and so on.
One particular subclass of non-Newtonian fluids is of Oldroyd type, that is,
$\mathbb{S}\triangleq \lambda\div u\,\Id+2\mu D(u)+\big( \f{W_{F}(F)F^T}{\det F}\big)$, where the deformation tensor $F$ satisfies the transport equation
\begin{equation}
\partial_tF+u\cdot\nabla F=\nabla uF.
\end{equation}
Formulations about viscoelastic flows of Oldroyd-B type are first introduced by Oldroyd \cite{O} and are extensively discussed in \cite{AM,La}. Consequently, we are concerned with the following compressible viscoelastic flow of Oldroyd type
\begin{equation}\label{R-E2}
\left\{\begin{array}{l}
\partial_t\rho+\div(\rho u)=0,\\[1ex]
\partial_t(\rho u)+\div(\rho u\otimes u)-\div\bigl(2\mu D(u)+\lambda\,\div u\, \Id\bigr)+\nabla\Pi=\textmd{div}\big( \f{W_{F}(F)F^T}{\det F}\big),\\[1ex]
\partial_tF+u\cdot\nabla F=\nabla uF.
\end{array}\right.
\end{equation}
where $W(F)$ is the elastic energy.  $ W_{F}(F)$ takes the Piola-Kirchhoff form and $\big( \f{W_{F}(F)F^T}{\det F}\big)$ is the Cauchy-Green tensor, respectively. For simplicity, a special form of the Hookean linear elasticity has been taken:
\be\label{he}
W(F)=\f{\al}{2}|F|^2,
\ee
where $ \al>0$ is elastic parameter.
The initial data are supplemented by
\be
(\rho,F;u)|_{t=0}=(\rho_0(x),F_0(x);u_0(x)), \q x\in\bR^n.  \label{id0}
\ee

In the present paper, we shall investigate the existence of global solution to the Cauchy problem \eqref{R-E2}-\eqref{id0}, as initial data are the perturbation of constant equilibrium state $(1,I,0)$. First of all, let us recall those previous efforts for incompressible viscoelastic flow, which reads as
\begin{equation}\label{R-E6}
\left\{\begin{array}{l}
\partial_t u+u\cdot\nabla u-\mu\Delta u+\nabla\Pi=\textmd{div}(FF^{T}),\\
\partial_tF+u\cdot\nabla F=\nabla uF, \\
\textmd{div} u=0.
\end{array}\right.
\end{equation}

For incompressible Oldroyd models, Renardy \cite{R} in 1985 investigated the existence and uniqueness of slow steady flows of viscoelastic fluids. The global existence of a small smooth solution was firstly established by Guillop\'{e} and Saut \cite{GS1}. Later, they \cite{GS2} investigated shearing motions and Poiseuille flows of Oldroyd fluids with retardation time, which exist for arbitrary time and arbitrary initial data. The case of $L^s$-$L^r$ solutions has been treated by Fernandez Cara, Guill\'{e}n and Ortega in \cite{FGO}. In higher dimensions, Lions and Masmoudi \cite{LM1} constructed global weak solutions for general initial conditions. Chemin and Masmoudi \cite{CM1} in the critical Besov space proved the existence and uniqueness of local and global solutions.
Constantin and Kliegel \cite{CK} established the global regularity of strong solutions for 2D Oldroyd-B fluids with additional diffusive stress.
Elgindi and Rousset \cite{ER} proved the global regularity of smooth solutions for 2D generalized Oldroyd-B type models without diffusive velocity. If the damping is absent in the classical Oldroyd case, the velocity viscosity alone may not be sufficient to guarantee the regularity of \eqref{R-E6}. The ``div-curl" structure is full explored by Lin, Liu and Zhang \cite{LLZ1}, Lei, Liu and Zhou \cite{LLZ2}, the Cauchy problem of \eqref{R-E6} admits the global classical solution in usual Sobolev spaces. Since then, there are a number of other results available in the assumption of structural conditions, see for example \cite{CZ1,LZ,ZF1,ZFZ}. Recently, the third author \cite{Zy1} in three dimensions proved the global existence of small solutions to
the incompressible Oldroyd-B model without damping mechanism. Her result can be also applied to the system \eqref{R-E6}, where the``div-curl" compatible condition is no longer needed. Recently, Chen and Hao \cite{CH} proved the global critical regularity in the Besov space based on the observation of Green's matrix. The reader is also referred to \cite{HLL,Lin} for the research summary of \eqref{R-E6}.

In this paper, we are concerned with the compressible viscoelastic flows. The mathematical modelling of compressible viscoelastic fluids was
proposed in earlier paper due to Beris and Edwards \cite{EB} (see also their book \cite{BE} or \cite{BP} and references therein). Fixed some positive time, Lei and Zhou \cite{LZ1} established the global existence of classical solutions to two-dimensional case, when initial data are subjected to incompressible constraints. Furthermore, the incompressible limit to \eqref{R-E6} was rigorously justified.
The existence and uniqueness of local-in-time strong solution with large initial data for the three-dimensional compressible viscoelastic flow was established by Hu and Wang \cite{HW4}. As the study of \eqref{R-E6}, the major difficulty proving the global existence of \eqref{R-E2} lies in the lacking of the dissipative estimates for the deformation and density. Inspired by the investigation of \eqref{R-E6} (see \cite{LLZ1,LLZ2}),  Hu-Wang \cite{HW1} and Qian-Zhang \cite{QZ1} independently explored intrinsic properties of \eqref{R-E2} such that the desired dissipation can be available. Indeed, their compatible conditions are listed as follows
\begin{equation}
 \rho_{0}\mathrm{det}F_{0}=1,\ \
\nabla\cdot(\rho_0 F_0^T)=0   \label{ip2}
\end{equation}
and
\begin{equation}
F_0^{lk}\partial_lF_0^{ij}-F_0^{lj}\partial_lF_0^{ik}=0. \label{ip22}
\end{equation}
The divergence constraint \eqref{ip2} makes sure that the gradient of $F$ behaves well in the elementary energy method, and \eqref{ip22} is used to control the quantity $\nabla\times F$. The conditional equivalence of  (\ref{ip2})-(\ref{ip22}) is shown by the recent work \cite{HZ}.

On the other hand, as in many works dedicated to compressible Navier-Stokes equatioins, \emph{scaling invariance}
plays a fundamental role. The reason why is that  whenever such an invariance exists,
  suitable critical quantities (that is, having the
same scaling invariance as the system under consideration) control the possible finite time blow-up, and
the global existence of strong solutions. Danchin \cite{Dr1} firstly solved (\ref{R-E1}) globally in critical homogeneous Besov spaces of $L^2$ type. Later, his result has been extended to those critical
Besov spaces that are not related to $L^2$, by Charve-Danchin \cite{CD} and Chen-Miao-Zhang \cite{CMZ} independently. Recently, Danchin and the second author \cite{DX,X} showed the optimal decay rates in general $L^p$ critical spaces. A natural (non Newtonian) extension in analysis is to consider (\ref{R-E2}).  Notice that (\ref{R-E2}) is invariant by the transformation
\begin{equation}\label{scaling}
\rho(t,x)\leadsto \rho(\ell ^2t,\ell x),\quad
u(t,x)\leadsto \ell u(\ell^2t,\ell x), \quad F(t,x) \leadsto F(\ell ^2t,\ell x) \qquad\ell>0,
\end{equation}
up to a change of the pressure term $\Pi$ into $\ell^2\Pi$ and the constant $\al$ into $\ell^2\al$. Under the assumptions \eqref{ip2}-\eqref{ip22}, Hu-Wang \cite{HW1} and Qian-Zhang \cite{QZ1} independently deduced \textit{a priori} dissipation estimates for complicated hyperbolic-parabolic systems, which lead to the existence of global existence in the critical $L^2$ Besov space. Hu-Wu \cite{HW2} proved the global existence of strong solutions to \eqref{R-E2} as initial data are the small perturbation $(1,I;0)$ in $H^2(\bR^3)$. Furthermore, it was shown that those solutions converged to equilibrium state at the decay rates of heat kernel. Barrett, Lu and E. S\"{u}li \cite{BLS} investigated 2D compressible Oldroyd-B type model which is derived from the compressible Navier-Stokes-Fokker-Planck system and proved the existence of large data global-in-time finite-energy weak solutions. Huo and Yong \cite{HY} studied the structural stability of a 1D compressible viscoelastic fluid model which was proposed by \"{O}ttinger \cite{O2} and established the global existence of smooth solutions near equilibrium.

Based on \cite{HW1,QZ1}, the first two authors \cite{PX1} established the global existence and time-decay estimates of solutions to \eqref{R-E2} in the general $L^p$ Besov space. The argument of effective velocities developed by Haspot \cite{Hb2}
was mainly employed, which is analogue of Hoff's viscous effective flux in \cite{Hd1}. Let us point out that the dissipation of
(\ref{R-E2}) with constraints \eqref{ip2}-\eqref{ip22} is standard, which is similar to that of the compressible Navier-Stokes equations (\ref{R-E1}). A question thus follows. \textit{Is it possible to find any new dissipative ingredients on non Newtonian compressible viscoelastic flows without \eqref{ip2}-\eqref{ip22}?} Here we aim at recasting the global-in-time existence of strong solutions in the framework of spatially Besov spaces with critical regularity without \eqref{ip2}-\eqref{ip22} that has been playing the key role in related efforts.

Before writing out the main statement of our paper, let us introduce some notation and definition first. To begin with, we need a Littlewood-Paley decomposition. There exists two radial smooth functions $\varphi(x),\chi(x)$ supported in the annulus $\mathcal{C}=\{\xi\in\bR^n:3/4\leq |\xi|\leq 8/3\}$ and the ball $B=\{\xi\in\bR^n:|\xi|\leq 4/3\}$, respectively such that
\be
\sum\limits_{j\in \bZ}\varphi(2^{-j}\xi)=1\q \forall \xi\in\bR^n\setminus\{0\}.\nn
\ee
The homogeneous dyadic blocks $\dD_j$ and the homogeneous low-frequency cut-off operators $\dot{S}_j$ are defined for all $j\in \bZ$ by
\be
\dD_j u=\varphi(2^{-j}D)f,\q \dot{S}_jf=\sum\limits_{k\leq j-1}\dD_k f=\chi(2^{-j}D)f.\nn
\ee
\q\ We denote by\ $\mathcal{Z}'(\bR^n)$ the dual space of
\be
\mathcal{Z}(\bR^n)\triangleq \{f\in \mathcal{S}(\bR^n):\p^\al \hat{f}(0)=0,\forall \al\in(\bN\cup 0)^n\}.\nn
\ee
Let us now turn to the definition of the main functional spaces and norms that
will come into play in our paper.
\begin{definition}
Let s be a real number and (p,r) be in $[1,\i]^2$. The homogeneous Besov space $\dot{B}^s_{p,r}$ consists of those distributions $u\in \mathcal{Z}'(\bR^n)$ such that
\be
\|u\|_{\dot{B}^s_{p,r}}\triangleq \Big(\sum\limits_{j\in\bZ}2^{jsr}\|\dD_j u\|^r_{L^p}\Big)^{\f{1}{r}}<\i.\nn
\ee
\end{definition}
Also, we introduce the hybrid Besov space since our analysis will be performed at different frequencies.
\begin{definition}
Let $s,\sigma\in \bR$. The hybrid Besov space $\B^{s,\sigma}$ is defined by
\be
\B^{s,\sigma}\triangleq \{f\in\mathcal{Z}'(\bR^n):\|f\|_{\B^{s,\sigma}}<\i\}, \nn
\ee
with
\be
\|f\|_{\B^{s,\sigma}}\triangleq \sum\limits_{2^k\leq R_0}2^{ks}\|\dD_k f \|_{L^2}+\sum\limits_{2^k>R_0}2^{k\sigma}\|\dD_k f \|_{L^2},\nn
\ee
\end{definition}
where $R_0$ is a fixed constant to be defined. $\B^{s,s}$ is the usual Besov space $\dB^s_{2,1}$ if $\sigma=s$. In the case where $u$ depends on the time variable, we consider the space-time mixed spaces as follows
\be
\|u\|_{L^q_T\B^{s,\sigma}}:=\big\|\|u(t,\cdot)\|_{\B^{s,\sigma}}\big\|_{L^q(0,T)}.\nn
\ee
In addition, we introduce another space-time mixed spaces, which is usually referred to Chemin-Lerner's spaces. The definition is given by
\be
\|u\|_{\t{L}^q_T\dot{B}^{s,\sigma}}\triangleq\sum\limits_{2^k\leq R_0}2^{ks}\|\dD_k u\|_{L^q(0,T)L^2}+\sum\limits_{2^k> R_0}2^{k\sigma}\|\dD_k u\|_{L^q(0,T)L^2}.\nn
\ee
The index $T$ will be omitted if $T=+\i$ and we shall denote by $\t{\mathcal{C}}_b(\B^{s,\sigma})$ the subset of functions $\t{L}^\i(\B^{s,\sigma})$ which are continuous from $\bR_+$ to $\B^{s,\sigma}$. It is easy to check that $\t{L}^1_T\B^{s,\sigma}=L^1_T\B^{s,\sigma}$ and $\t{L}^q_T\B^{s,\sigma}\subseteq L^q_T\B^{s,\sigma}$ for $q>1$.

Our results are stated as follows.

\begin{theorem}\label{thglobal1}
Let $I$ be the unit matrix of order $n$.
 There exists two positive constants $\eta$ and $M$ such that if
\be
(\rho_0-1,F_0-I;u_0)\in\Big( \B^{n/2-1,n/2}\Big)^{1+n^2}\times \Big(\dB^{n/2-1}_{2,1}\Big)^n. \nn
\ee
and
\be
\|(\rho_0-1,F_0-I)\|_{\B^{n/2-1,n/2}}+\|u_0\|_{\dB^{n/2-1}_{2,1}}\leq \eta,\nn
\ee
then the Cauchy problem \eqref{R-E2} and \eqref{id0} has a global unique solution $(\rho,F;u)$ such that
\bea
&&u\in\Big(\t{\mathcal{C}_b}(\bR_+;\dB^{n/2-1}_{2,1})\cap L^1(\bR_+;\dB^{n/2+1}_{2,1})\Big)^{n}\nn\\
&&(\rho-1)\in \t{\mathcal{C}_b}(\bR_+;\B^{n/2-1,n/2});(F-I)\in \Big(\t{\mathcal{C}_b}(\bR_+;\dB^{n/2}_{2,1})\Big)^{n^2}.
\eea
Moreover, the following estimate holds:
\bea\label{1.4}
\|\rho-1\|_{\t{L}^\i\B^{n/2-1,n/2}}+\|F-I\|_{\t{L}^\i\dB^{n/2}_{2,1}}+\|u\|_{\t{L}^\i\dB^{n/2-1}_{2,1}\cap{L}^1\dB^{n/2+1}_{2,1}} \leq M\eta.
\eea
\end{theorem}
\begin{remark}
We would like to mention that the third author \cite{Zy2} got a global-in-time existence of smooth solutions in Sobolev space without
(\ref{ip2})-(\ref{ip22}), from which one can see that there is some regularity loss on both density and deformation tensor. In the Besov framework,
we prove the evolution of critical regularity of perturbation variable $(\rho, u, \tau)$ with $\tau\triangleq\f{FF^T}{\det F}-I$:
\bea\label{1.5}
&&\|(\rho-1, \tau)\|_{\t{L}^\i\B^{n/2-1,n/2}}+\|u\|_{\t{L}^\i\dB^{n/2-1}_{2,1}\cap{L}^1\dB^{n/2+1}_{2,1}} \leq M\eta.
\eea
Therefore, our analysis allows to establish the corresponding global-in-time result for the compressible Oldroyd-B model without damping. Back to the original
system \eqref{R-E2}, one resorts to the estimate of transport equation, which leads to one regularity loss of $F$ at low frequencies (see \eqref{4.56} for details).
\end{remark}
\begin{remark}
In absence of (\ref{ip2})-(\ref{ip22}), the new effective flux $(\theta,\mathcal{G})$ (see below) plays the key role in the analysis, which enables us to capture the partial dissipation arising from the combination of density and deformation only (see (\ref{1.4}) or (\ref{1.5})).
Consequently, the density and deformation tensor themselves might grow in time, which is totally different in comparison with those efforts for compressible Navier-Stokes equations \eqref{R-E1} (see for example \cite{Dr1,CD,CMZ,Hb2}).
\end{remark}

\begin{remark}
It is possible to prove the analogue of Theorem \ref{thglobal1} in more general $L^p$ framework. This is beyond our primary interest in the present paper, since we focus on the elementary dissipative structure of non-Newtonian fluids. Actually, the $L^2$ orthogonal property of projection operator $(\mP,\, \mvP)$ is well used in the proof of Proposition \ref{proapriori}, see \eqref{4.36} and \eqref{4.59} for details.
\end{remark}

We end this section with a strategy in the proof of Theorem \ref{thglobal1}. The starting point is to
rewrite \eqref{R-E2} as the linearized compressible viscoelastic flows about $(1,I,0)$. In order to avoid those initial compatible conditions, one can
view $\f{FF^T}{\det F}$ as a new variable rather than the nonlinear term in previous efforts. Without loss of generality, we set $P'(1)=1$. Define
\be
\rho=1+a,\q  p(t,x)=\Pi(1+a)-\Pi(1),\q \tau(t,x)=\f{FF^T}{\det F}(t,x)-I.\nn
\ee
It is shown that by the direct computation
 \be\label{r1}
\left\{
\begin{aligned}
&\p_t p+u\cdot\nabla p+\nabla\cdot u=F_{1}, \\
&\p_tu-\mathcal{A}u+(\na p-\al\na\cdot\tau)=F_{2},\\
&\p_t\tau+u\cdot \na \tau+\na\cdot u \text{Id}-2D(u)=F_{3},
\end{aligned}
\right.
\ee
with \be F_{1}\triangleq -K(a)\nabla\cdot u \nn,
\ee
\be
F_2\triangleq-I(a)\mathcal{A}u-u\cdot\nabla u+I(a)(\na p-\al\na\cdot\tau)+\f{1}{1+a}\text{div}\big(2\t{\mu}(a)D(u)+{\la}(a)\text{div} u\text{Id}\big)\nn\\
\ee
and
\be
F_3 \triangleq \na u\tau+\tau(\na u)^T-\na\cdot u\tau,
\nn
\ee
where
\be
I(a)\triangleq\f{a}{1+a},\ K(a)\triangleq \Pi'(1+a)(1+a)-1,\ \mathcal{A}=\mu(1)\Dl+(\la(1)+\mu(1))\nabla\text{div},
\nn
\ee
and
\be\ \t{\mu}(a)=\mu(1+a)-\mu(1),\  \t{\la}(a)=\la(1+a)-\la(1).\nn
\ee
For simplicity, we denote $\la(1)=\la_0, \mu(1)=\mu_0$. In order to capture the dissipation arising from the complicated coupling between $p$ and $\tau$, let us introduce the \textit{new effective flux} $\theta=\na p-\al\na\cdot\tau$. By employing the operator $\nabla$ to $\eqref{r1}_1$ and the operator $\alpha\mathrm{div}$ to $\eqref{r1}_3$, respectively,  one can get
\be\label{r2}
\left\{
\begin{aligned}
&\p_t \th+u\cdot\na\th+\al\Dl u+\na(\nabla\cdot u)=\t{F}_1-\al\t{F}_3, \\
&\p_t u-\mathcal{A} u+\th= F_2\\
\end{aligned}
\right.
\ee
with $\t{F}_1=-\na(K(a)\na\cdot u)-(\na u)^T\na p,\t{F}_3=-((\na u)^T\cdot\na)\cdot \tau+\na\cdot F_3$. The corresponding linear system reads
\be\label{line}
\left\{
\begin{aligned}
&\p_t \th+\al\Dl u+\na(\nabla\cdot u)=0, \\
&\p_t u-\mathcal{A} u+\th=0,\\
\end{aligned}
\right.
\ee

As shown by the formal spectral analysis in Section 2, we see that \eqref{line} admits the similar dissipation structure as that of usual compressible Navier-Stokes equations. The observation on the combination of density and deformation tensor (\textit{without any compatible conditions})  is new in compressible non-Newtonian fluids, which enables us to establish a global existence in the critical Besov space.

\section{Formal spectral analysis and energy functionals}
In order to understand the proof of Theorem \ref{thglobal1}, it is convenient to give the formal spectral analysis for \eqref{line}.
For $s\in \bR$, we denote $\Lambda^sf\triangleq\mathcal{F}^{-1}(|\xi|^s\mathcal{F}(f))$. Also we use $\mP$ to denote the projection operator $I+(-\Dl)^{-1}\na\na\cdot$ and $\mvP=-(-\Dl)^{-1}\na\na\cdot$  on the divergence-free vector and potential vector, respectively. By applying $\La^{-1}\mP,\La^{-1}\mvP$ to the first equation of \eqref{line} and $\mP,\mvP$ to the second equation of \eqref{line}, we get
\be\label{elinear}
\left\{
\begin{aligned}
&\p_t \La^{-1}\mvP\th -(1+\alpha)\La \mvP u= 0, \\
&\p_t\mvP u-(\lambda_{0}+2\mu_{0})\Dl\mvP u+\mvP\th=0,\\
&\p_t\La^{-1}\mP\th-\alpha\La \mP u=0,\\
&\p_t\mP u-\mu_{0}\Dl\mP u+\mP\th=0.\\
\end{aligned}
\right.
\ee
Clearly, we see that there are \textit{two hyperbolic-parabolic coupled systems} for $(\La^{-1}\mvP\th,\mvP u)$ and $(\La^{-1}\mP\th, \mP u)$ available, which are similar with the case of compressible Navier-Stokes equations (see \cite{Dr1}). For example,
we investigate the $2\times 2$ subsystem for $(\Phi,\Psi)\triangleq (\La^{-1}\mP\th, \mP u)$:
\be\label{elinear1}
\left\{
\begin{aligned}
&\p_t\Phi-\alpha\La \Psi=0,\\
&\p_t\Psi-\mu_{0}\Dl \Psi+\Lambda \Phi=0.\\
\end{aligned}
\right.
\ee
The Green matrix is given by $G(D)=\left(
                         \begin{array}{cc}
                           0 & \alpha \Lambda \\
                           -\Lambda & \mu_{0}\Delta \\
                         \end{array}
                       \right).
$
Let $\la_\pm$ be the eigenvalues of $G(\xi)$. For low frequencies $(\mu_{0}|\xi|<2\sqrt{\alpha})$, the eigenvalues are
\[
\la_\pm=-\frac{\mu_{0}|\xi|^2}{2}\Big(1\pm\mathrm{i}\sqrt{\frac{4\alpha}{\mu_{0}^2|\xi|^2}-1}\Big).
\]
The situation of high frequencies $(\mu_{0}|\xi|>2\sqrt{\alpha})$ is quite different. The eigenvalues are
\[
\la_\pm=-\frac{\mu_{0}|\xi|^2}{2}\Big(1\pm\sqrt{1-\frac{4\alpha}{\mu_{0}^2|\xi|^2}}\Big).
\]
Consequently, as in \cite{Dr1}, one can expect a parabolic smoothing for low frequencies of $(\Phi,\Psi)$, a damping for high frequencies of $\Phi$ and a parabolic smoothing
for high frequencies of $\Psi$. Similar analysis can be performed for another hyperbolic-parabolic system. So it is reasonable to define the energy at low frequencies as
\be
\t{\mE}^\ell_T:=\sup\limits_{t\in[0,T)}\|(u,\La^{-1}\th)\|^\ell_{\dB^{n/2-1}_{2,1}}+\int^T_0\|(u,\La^{-1}\th)\|^\ell_{\dB^{n/2+1}_{2,1}}dt,
\ee
and the energy at high frequencies
\be
\t{\mE}^h_T:=\sup\limits_{t\in[0,T)}\lt(\|u\|^h_{\dB^{n/2-1}_{2,1}}+\|\La^{-1}\th\|^h_{\dB^{n/2}_{2,1}}\rt)+\int^T_0\|u\|^h_{\dB^{n/2+1}_{2,1}}dt+\int^T_0\|\La^{-1}\th\|^h_{\dB^{n/2}_{2,1}}dt.
\ee

The above analysis looks so standard, however, keep in mind that the partial dissipation of density and deformation tensor is captured only. As a matter of fact, we have to meet those nonlinear terms (see $F_1,\ F_2$ and $F_3$) with respect to the variable $(a,\tau)$ itself. In order to close the energy method,  we need additional $L^\i$ estimates for $a$ and $\tau$ in time. For that end, we introduce another \textit{new effective flux} $\mathcal{G}\triangleq\tau-p \text{Id}$. It follows from
$\eqref{r1}_1 \text{Id}$ and $\eqref{r1}_3$ that
 \be\label{r4}
\left\{
\begin{aligned}
&\p_t p+u\cdot\na p+\nabla\cdot u=F_1, \\
&\p_tu-\mathcal{A}u+(1-\al)\na p-\al\na\cdot\mathcal{G}=F_2,\\
&\p_t\mG+u\cdot\na \mathcal{G}-2D(u)=F_3-F_1\text{Id}.
\end{aligned}
\right.
\ee
Furthermore, we revise our energy functionals (2.3)-(2,4) a little bit, which are given by
\be\label{el}
{\mE}^\ell_T:=\sup\limits_{t\in[0,T)}\|(p,\tau,\ u,\La^{-1}\th)\|^\ell_{\dB^{n/2-1}_{2,1}}+\int^T_0\|(u,\La^{-1}\th)\|^\ell_{\dB^{n/2+1}_{2,1}}dt,
\ee
and
\bea\label{eh}
\hspace{-5mm} {\mE}^h_T&:=&\sup\limits_{t\in[0,T)}\lt(\|u\|^h_{\dB^{n/2-1}_{2,1}}+\|(p,\ \tau,\ \La^{-1}\th)\|^h_{\dB^{n/2}_{2,1}}\rt)\nn\\
&&\qq\qq+\int^T_0\|u\|^h_{\dB^{n/2+1}_{2,1}}dt+\int^T_0\|\La^{-1}\th\|^h_{\dB^{n/2}_{2,1}}dt.
\eea
Indeed, by combining the $L^\i$ estimates of $(p,\tau)$ (see \eqref{4.15}, \eqref{e4.4} and \eqref{e4.18} for details), it is sufficient to establish the global existence of strong solutions to the Cauchy problem \eqref{R-E2}-\eqref{id0}. Finally, it's worth noting that our analysis holds true for non-small coupling parameter $\alpha$.

\section{A priori energy estimates}

Following from the spectral analysis in Section 2, we shall prove crucial \textit{a priori} estimates for the energy functionals in \eqref{el} and \eqref{eh}.

Let $T>0$. We by $\mathcal{E}_T$ denote the functional space
\bea
\mathcal{E}_T&\triangleq& \big\{(p,\tau;u,\th)|(p,\tau)\in\t{L}^\i(0,T;\B^{n/2-1,n/2}),\nn\\
             &&\q u\in\t{L}^\i(0,T;\dB^{n/2-1}_{2,1})\cap L^1(0,T;\dB^{n/2+1}_{2,1}), \nn\\
             &&\q \La^{-1}\th\in\t{L}^\i(0,T;\B^{n/2-1,n/2})\cap L^1(0,T;\B^{n/2+1,n/2})\big\}\nn
\eea
and the corresponding norm is given by
\bea
&&\|(p,\tau;u,\th)\|_{\mathcal{E}_T}\triangleq \|(p,\tau)\|_{\t{L}^\i_T\B^{n/2-1,n/2}}\nn\\
&&\qq\qq+\|u\|_{\t{L}^\i_T\dB^{n/2-1}_{2,1}\cap L^1_T\dB^{n/2+1}_{2,1}}+\|\La^{-1}\th\|_{\t{L}^\i_T\B^{n/2-1,n/2}\cap L^1_T\B^{n/2+1,n/2}}.\nn
 \eea
Note that $p=\Pi(1+a)-\Pi(1),\ \Pi'(1)=1$. There exists a small number $\eta_0>0$ such that $\|a\|_{L^\i([0,T]\times\bR^n)}\leq \eta_0$. Consequently,
 $a$ can be expressed by a smooth function of $p$. Set $a=h(p)$.
\begin{proposition}\label{proapriori}
 Assume that $(p,\tau;u)$ is a strong solution of System \eqref{r1} on $[0,T]$ with
\be
\|a\|_{L^\i([0,T]\times\bR^n)}\leq \eta_0.  \nn
\ee
Then it holds that
\bea\label{4.1}
&&\|(p,\tau;u,\th)\|_{\E_T}\leq C\Big\{\|(p,\tau;u)(0)\|_{\E_0}\nn\\
&&\qq\qq\qq\qq+\|(p,\tau;u,\th)\|^2_{\E_T}\big(1+\|(p,\tau;u,\th)\|_{\E_T}\big)^{n+3}\Big\},
\eea
where $\|(p,\tau;u)(0)\|_{\E_0}\triangleq \|(p,\tau)(0)\|_{\B^{n/2-1,n/2}}+\|u(0)\|_{\dB^{n/2-1}_{2,1}}$.
\end{proposition}

We divide the proof of Proposition  \ref{proapriori} into three parts for clarity. The first part is devoted to dissipative estimates for variable $(\th,u)$. More precisely, the parabolic smoothing effect for low frequencies of $(\th,u)$, the damping for high frequencies of $\th$ and the parabolic smoothing effect for high frequencies of $u$ will be addressed. With the help of the new effective flux $\mathcal{G}$, in the second part, we give the additional
$L^\i$ estimates for full variables $p$ and $\tau$ in time. The last part is dedicated to bounding of those nonlinear terms.

\subsection{Dissipative estimates of $\boldsymbol{(\th,u)}$ }

\q In this subsection, we derive the parabolic smoothing effect for low frequencies of $(\th,u)$, the damping for high frequencies of $\th$ and the parabolic smoothing effect for high frequencies of $u$.

Set $(\Phi,\Psi)\triangleq (\La^{-1}\mP\th, \mP u)$ and $(\Phi^{\bot},\Psi^{\bot})\triangleq (\La^{-1}\mvP\th, \mvP u)$. We use the notation
$f_k=\dD_k f$ for any scalar (vector or matrix, respectively) function $f$. \\

\noindent\textbf{Step 1: Low-frequency estimates ($\boldsymbol{2^k\leq R_0}$) }\\
By applying $\La^{-1}\mP\dD_k,\La^{-1}\mvP\dD_k$ to the first equation of \eqref{r2}, $\mP\dD_k,\mvP\dD_k$ to the second equation of \eqref{r2}, we have
\be\label{r6}
\left\{
\begin{aligned}
&\p_t \Phi^\bot_k+u\cdot\na\Phi^\bot_k -(1+\al)\La \Psi^\bot_k= \La^{-1}\mvP\dD_k(\t{F}_1-\al\t{F}_3)+\mathcal{R}^1_k, \\
&\p_t\Psi^\bot_k-(\la_0+2\mu_0)\Dl\Psi^\bot_k+\La\Phi^\bot_k=\dD_k\mvP F_2,\\
&\p_t\Psi_k-\mu_0\Dl\Psi_k+\La\Phi_k=\dD_k\mP F_2,\\
&\p_t\Phi_k+u\cdot\na\Phi_k-\al\La \Psi_k=\La^{-1}\mP\dD_k(\t{F}_1-\al\t{F}_3)+\mathcal{R}^2_k,
\end{aligned}
\right.
\ee
 where commutators are given by $\mathcal{R}^1_k=[u\cdot\na,\La^{-1}\mvP\dD_k]\th,\ \mathcal{R}^2_k=[u\cdot\na,\La^{-1}\mP\dD_k]\th$.

 Taking $L^2$ inner product of $\eqref{r6}_1$  with $\f{1}{1+\al}\Phi^\bot_k$, $\eqref{r6}_2$ with $\Psi^\bot_k$, $\eqref{r6}_3$ with $\Psi_k$ and $\eqref{r6}_4$ with $\f{1}{\al}\Phi_k$ respectively, and then adding  the resulting equations together, we obtain
\bea\label{4.5}
&&\f{1}{2}\f{d}{dt}\Big(\| (\Psi_k,\Psi^\bot_k)\|^2_{L^2}+\f{1}{\al}\|\Phi_k\|^2_{L^2}+\f{1}{1+\al}\|\Phi^\bot_k\|^2_{L^2}\Big)\nn\\
&& \hspace{5mm} +\mu_0\|\La\Psi_k\|^2_{L^2}+(\la_0+2\mu_0)\|\La\Psi^\bot_k\|^2_{L^2}\nn\\
&&=( u_k|\dD_k F_2)+\big(\f{1}{2(1+\al)}|\Phi^\bot_k|^2+\f{1}{2\al}|\Phi_k|^2\big|\na\cdot u\big)+\f{1}{\al}\big(\Phi_k|\La^{-1}\mP \dD_k(\t{F}_1-\al\t{F}_3)+\mathcal{R}^2_k\big)\nn\\
&&\q +\f{1}{1+\al}\big(\Phi^\bot_k|\La^{-1}\mvP\dD_k(\t{F}_1-\al\t{F}_3)+\mathcal{R}^1_k\big).
\eea

To capture the dissipation arising from $\th$, we take the $L^2$ inner product of $\eqref{r6}_1$ with $-\La\Psi^\bot_k$, $\eqref{r6}_2$ with $\La\Phi^\bot_k$, $\eqref{r6}_3$ with $\La\Phi_k$ and $\eqref{r6}_4$ with $-\La\Psi_k$ respectively. Then we add these resulting equations together and get
\bea\label{4.6}
 &&\f{d}{dt}\big[(\Psi^\bot_k|\La\Phi^\bot_k)+(\Psi_k|\La\Phi_k)\big]+\|\La(\Phi_k,\Phi^\bot_k)\|^2_{L^2}+\al\|\La \Psi_k\|^2_{L^2}+(1+\al)\|\La \Psi^\bot_k\|^2_{L^2}\nn\\
 &&\hspace{5mm}-\mu_0(\Dl\Psi_k|\La\Phi_k)-(\la_0+2\mu_0)(\Dl\Psi^\bot_k|\La\Phi^\bot_k)\nn\\
&=&-\big(\La\Psi^\bot_k|\La^{-1}\mvP\dD_k(\t{F}_1-\al\t{F}_3)+\mathcal{R}^1_k\big)+(\th_k|\dD_k F_2)\nn\\
&&-\big(\La\Psi_k|\La^{-1}\mP\dD_k(\t{F}_1-\al\t{F}_3)+\mathcal{R}^2_k\big)+(u\cdot\na\La^{-1}\th_k\big|\La u_k).
\eea
Now, we multiply a small constant $\nu$ (to be determined) to \eqref{4.6} and then add
the resulting equation with \eqref{4.5} together. Consequently, we are led to the following inequality
\bea\label{e4.9}
 &&\f{d}{dt}f^2_{\ell,k}+\t{f}^2_{\ell,k}\ls |F(t)|,
\eea
where \bea
&&f^2_{\ell,k}:=\|(\Psi_k,\Psi^\bot_k)\|^2_{L^2}+\f{1}{\al}\|\Phi_k\|^2_{L^2}+\f{1}{1+\al}\|\Phi^\bot_k\|^2_{L^2}\nn\\
 &&\qq\qq+2\nu(\Psi^\bot_k|\La\Phi^\bot_k)+2\nu(\Psi_k|\La\Phi_k),\nn\\
&&\t{f}^2_{\ell,k}:=(\mu_0+\al\nu)\|\La \Psi_k\|^2_{L^2}+(\la_0+2\mu_0+(1+\al)\nu)\|\La \Psi^\bot_k\|^2_{L^2}+\nu\|\La(\Phi_k,\Phi^\bot_k)\|^2_{L^2}\nn\\
&&\qq-\mu_0\nu(\Dl\mP u_k|\mP\th_k)-(\la_0+2\mu_0)\nu(\Dl\mvP u_k|\mvP\th_k).\nn
\eea
and
\bea
&&F(t):=
( u_k|\dD_k F_2)+\big(\f{1}{2(1+\al)}|\Phi^\bot_k|^2+\f{1}{2\al}|\Phi_k|^2\big|\na\cdot u\big)
+\f{1}{\al}\big(\Phi_k|\La^{-1}\mP \dD_k(\t{F}_1-\al\t{F}_3)\nn\\
 &&\qq+\mathcal{R}^2_k\big)
+\f{1}{1+\al}\big(\Phi^\bot_k|\La^{-1}\mvP\dD_k(\t{F}_1-\al\t{F}_3)+\mathcal{R}^1_k\big)
-\nu\big(\La\Psi^\bot_k|\La^{-1}\mvP\dD_k(\t{F}_1-\al\t{F}_3)\nn\\
 &&\qq+\nu\mathcal{R}^1_k\big)+(\th_k|\dD_k F_2)
-\nu\big(\La\Psi_k|\La^{-1}\mP\dD_k(\t{F}_1-\al\t{F}_3)+\nu\mathcal{R}^2_k\big)+(u\cdot\na\La^{-1}\th_k\big|\La u_k).\nn
\eea
For any fixed $R_0$, we choose  $\nu\sim\nu(\la_0,\mu_0,R_0)$ sufficiently small such that
\be
\begin{aligned}
&f^2_{\ell,k}\sim \|u_k\|^2_{L^2}+\|\La^{-1}\th_k\|^2_{L^2},\\
&\t{f}^2_{\ell,k}\sim 2^{2k}\big(\|u_k\|^2_{L^2}+\|\La^{-1}\th_k\|^2_{L^2}\big).
\end{aligned}
\ee
By using Cauchy-Schwarz inequality, furthermore, one can get owing to $2^k\leq R_0$,
\bea\label{4.25}
&&\f{d}{dt}f_{\ell,k}+2^{2k}f_{\ell,k}\ls \|\dD_k (\La^{-1}\t{F}_1,\La^{-1}\t{F}_3, F_2),\mathcal{R}^1_k,\mathcal{R}^2_k\|_{L^2}\nn\\
&&\qq\qq\qq\qq  +\|\La^{-1}\th_k\na\cdot u,u\cdot\na\La u_k,\na\cdot u\La u_k\|_{L^2}, \eea
which indicates that
\bea\label{4.11}
&&\q\|( u,\La^{-1}\th)\|^{\ell}_{{\t{L}^\i_T\dB^{\frac{n}{2}-1}_{2,1}}\cap \t{L}^1_T\dB^{\frac{n}{2}+1}_{2,1}}\nn\\
&&\ls \|(u,\La^{-1}\th)(0)\|^{\ell}_{\dB^{\frac{n}{2}-1}_{2,1}}
+\|(\La^{-1}\t{F}_1, F_2,\La^{-1}\t{F}_3)\|^\ell_{\t{L}^1_T\dB^{\frac{n}{2}-1}_{2,1}}\nn\\
&&\hspace{5mm} +\sum\limits_{2^j\leq R_0}2^{j(n/2-1)}\int^T_0\Big\|\mathcal{R}^1_k,\mathcal{R}^1_k,\La^{-1}\th_k\na\cdot u,u\cdot\na\La u_k
,\na\cdot u\La u_k\Big\|_{L^2}dt\nn\\
&&\ls \|(p,\tau;u)(0)\|^{\ell}_{\dB^{\frac{n}{2}-1}_{2,1}}
+\|(\La^{-1}\t{F}_1, F_2,\La^{-1}\t{F}_3)\|^\ell_{\t{L}^1_T\dB^{\frac{n}{2}-1}_{2,1}}\nn\\
&&\hspace{5mm}+\sum\limits_{2^j\leq R_0}2^{j(n/2-1)}\int^T_0\Big\|\mathcal{R}^1_k,\mathcal{R}^2_k,\La^{-1}\th_k\na\cdot u,u\cdot\na\La u_k
,\na\cdot u\La u_k\Big\|_{L^2}dt.
\eea
It follows from those commutator estimates in \cite{BCD1} that
\be
\sum\limits_{k\in\bZ}2^{ks}\int^T_0\|\mathcal{R}^i_k\|_{L^2}dt\ls \|\na u\|_{\t{L}^{r_1}_T\dB^{n/2}_{2,1}}
\|\La^{-1}\th\|_{\t{L}^{r_2}_T\dB^{s}_{2,1}}. \q \text{for}\q i=1,2.\nn
\ee
where $\f{1}{r_1}+\f{1}{r_2}=1.$
In particular, one can get
\bea\label{e4.10e}
\sum\limits_{2^k\leq R_0}2^{k(n/2-1)}\int^T_0\big\|\mathcal{R}^1_k,\mathcal{R}^2_k\big\|_{L^2}dt\ls \|\na u\|_{\t{L}^{1}_T\dB^{n/2}_{2,1}}\|\La^{-1}\th\|_{\t{L}^{\i}_T\dB^{n/2}_{2,1}}\ls\|(p,\tau;u,\th)\|^2_{\mathcal{E}_T}.
\eea
Similarly, regarding other terms in the last integral of \eqref{4.11}, we arrive at
\bea\label{e4.11e}
\hspace{-10mm}\sum\limits_{2^k\leq R_0}2^{k(n/2-1)}\int^T_0\big\|\La^{-1}\th_k\na\cdot u\big\|_{L^2}dt
&\ls& \int^T_0\|\na u\|_{L^\i}dt\|\La^{-1}\th\|^\ell_{\t{L}^{\i}_T\dB^{n/2-1}_{2,1}}\nn\\
&\ls &\|\na u\|_{\t{L}^{1}_T\dB^{n/2}_{2,1}}\|\La^{-1}\th\|_{\t{L}^{\i}_T\dB^{n/2}_{2,1}}\nn\\
&\ls&\|(p,\tau;u,\th)\|^2_{\mathcal{E}_T},
\eea
\bea\label{e4.12e}
\sum\limits_{2^k\leq R_0}2^{k(n/2-1)}\int^T_0\big\|u\cdot\na\La u_k\big\|_{L^2}dt
&\ls &\big(\int^T_0\|u\|^2_{L^\i}dt\big)^{1/2}\|\La^{2}u_k\|^\ell_{\t{L}^{2}_T\dB^{n/2-1}_{2,1}}\nn\\
&\ls& \| u\|_{\t{L}^2_T\dB^{n/2}_{2,1}}\|u\|^\ell_{\t{L}^{2}_T\dB^{n/2}_{2,1}}\nn\\
&\ls&\|(p,\tau;u,\th)\|^2_{\mathcal{E}_T}
\eea
and
\bea\label{e4.13e}
\hspace{-10mm}\sum\limits_{2^k\leq R_0}2^{k(n/2-1)}\int^T_0\big\|\na\cdot u\La u_k\big\|_{L^2}dt
&\ls& \int^T_0\|\na u\|_{L^\i}dt\|\La u\|^\ell_{\t{L}^{\i}_T\dB^{n/2-1}_{2,1}}\nn\\
&\ls& \| u\|_{\t{L}^1_T\dB^{n/2+1}_{2,1}}\|u\|^\ell_{\t{L}^{\i}_T\dB^{n/2-1}_{2,1}}\nn\\
&\ls&\|(p,\tau;u,\th)\|^2_{\mathcal{E}_T}.
\eea
Together with \eqref{4.11}-\eqref{e4.13e}, we deduce that
\bea\label{e4.14}
\|( u,\La^{-1}\th)\|^{\ell}_{{\t{L}^\i_T\dB^{\frac{n}{2}-1}_{2,1}}\cap \t{L}^1_T\dB^{\frac{n}{2}+1}_{2,1}}
&\ls& \|(p,\tau;u)(0)\|^{\ell}_{\dB^{\frac{n}{2}-1}_{2,1}} +\|(p,\tau;u,\th)\|^2_{\mathcal{E}_T}\nn\\
&&+\|(\La^{-1}\t{F}_1, F_2,\La^{-1}\t{F}_3)\|^\ell_{\t{L}^1_T\dB^{\frac{n}{2}-1}_{2,1}}. \eea

\noindent\textbf{Step 2: High-frequency estimates ($\boldsymbol{2^k> R_0}$) }\\
At high frequencies, let us perform the effective velocity argument as in \cite{Hb2} that was originated from
Hoff's viscous effective flux in \cite{Hd1}, and overcome the loss of one derivative of $(\Phi,\Phi^{\bot}) $. By applying $\mP$ and $\mvP$ to \eqref{r2}, we get
\be\label{r9}
  \left\{
\begin{aligned}
&\p_t\mP\th+\al\Dl \mP u=\mP(\t{F}_1-\al\t{F}_3-u\cdot\na\th\big), \\
&\p_t\mP u-\mu_0 \Dl \mP u+\mP\th=\mP F_2,\\
&\p_t\mvP\th+(1+\al)\Dl \mvP u=\mvP(\t{F}_1-\al\t{F}_3-u\cdot\na\th\big), \\
&\p_t\mvP u-(\la_0+2\mu_0) \Dl \mvP u+\mvP\th=\mvP F_2.
\end{aligned}
\right.
\ee
We define the effective velocities $w, {w}^\bot$ such that  $-\mu_0 \Dl \mP u+\mP\th=-\mu_0\Dl w$ and $-(\la_0+2\mu_0) \Dl \mvP u+\mvP\th=-(\la_0+2\mu_0)\Dl {w}^\bot$. It follows that
\be\label{e4.48}
 w=\mP u+\f{1}{\mu_0}(-\Dl)^{-1} \mP\th,\ {w}^\bot=\mvP u+\f{1}{\la_0+2\mu_0}(-\Dl)^{-1} \mvP\th.
\ee
Firstly, we do some estimates for effective velocities $w$ and ${w}^\bot$. It is easy to check that
\be\label{4.50}
  \left\{
\begin{aligned}
&\p_t w-\mu_0 \Dl w+\f{\al}{\mu^2_0}(-\Dl)^{-1}\mP \th-\f{\al}{\mu_0}w\\
&=\mP F_2+\f{1}{\mu_0}(-\Dl)^{-1}\mP\big(\t{F}_1-\al\t{F}_3-u\cdot\na\th\big),\\
&\p_t{w}^\bot-(\la_0+2\mu_0) \Dl{w}^\bot+\f{1+\al}{(\la_0+2\mu_0)^2}(-\Dl)^{-1}\mvP\th-\f{1+\al}{\la_0+2\mu_0}{w}^\bot\\
&=\mvP F_2+\f{1}{\la_0+2\mu_0}(-\Dl)^{-1}\mvP\big(\t{F}_1-\al\t{F}_3-u\cdot\na\th\big).
\end{aligned}
\right.
 \ee
Note that \eqref{2.3} to \eqref{4.50}, the parabolic smooth estimate (See Lemma \ref{lheat}) enables us to obtain
\be\label{4.51}
\begin{aligned}
&\|(w,{w}^\bot)\|^h_{\t{L}^\i_T\dB^{n/2-1}_{2,1}\cap L^1_T\dB^{n/2+1}_{2,1}}\ls \|(w_0,{w}^\bot_0)\|^h_{\dB^{n/2-1}_{2,1}}
+\|(w,{w}^\bot)\|^h_{L^1_T\dB^{n/2-1}_{2,1}}+\|
\th\|^h_{L^1_T\dB^{n/2-3}_{2,1}}\\
&\q\qq+\|(\t{F}_1,\t{F}_3)\|^h_{L^1(\dB^{n/p-3}_{p,1})}+\|F_2\|^h_{L^1(\dB^{n/2-1}_{2,1})}+\|u\cdot\na\th\|^h_{L^1(\dB^{n/2-3}_{2,1})}.\\
\end{aligned}
\ee
Owing to the high frequency cut-off $2^k>R_0$, we have
\be
\|(w,{w}^\bot)\|^h_{L^1_T\dB^{n/2-1}_{2,1}}\ls R^{-2}_0\|(w,{w}^\bot)\|^h_{L^1_T\dB^{n/2+1}_{2,1}}, \q \|\th\|^h_{L^1_T\dB^{n/p-3}_{p,1}}\ls
R^{-2}_0\|\th\|^h_{L^1_T\dB^{n/2-1}_{1,1}}.\nn
\ee

Choosing $R_0>0$ sufficiently large, the terms $\|(w,{w}^\bot)\|^h_{L^1_T\dB^{n/2-1}_{2,1}}$ on right-side of \eqref{4.51}
 can be absorbed by the corresponding parts on left-hand side of \eqref{4.51}. Consequently, we conclude that
\be\label{4.52}
\begin{aligned}
&\|(w,{w}^\bot)\|^h_{\t{L}^\i_T\dB^{n/2-1}_{2,1}\cap L^1_T\dB^{n/2+1}_{2,1}}\ls \|(w_0,{w}^\bot_0)\|^h_{\dB^{n/2-1}_{2,1}}
+R^{-2}_0\|\th\|^h_{L^1_T\dB^{n/2-1}_{2,1}}\\
&\q\qq+R^{-2}_0\|(\t{F}_1,\t{F}_3)\|^h_{L^1(\dB^{n/2-1}_{2,1})}+\|F_2\|^h_{L^1(\dB^{n/2-1}_{2,1})}+\|u\cdot\na\th\|^h_{L^1(\dB^{n/2-3}_{2,1})}.\\
\end{aligned}
\ee

Next, we intend to obtain the damping estimate for $\theta$ at high frequencies. Indeed, it follows from the first equation of \eqref{r2} that
\be\label{4.39}
\begin{aligned}
&\p_t\th+u\cdot\na \th +\al\Dl u+\na\na\cdot u=\t{F}_1-\al\t{F}_3,
\end{aligned}
\ee
Applying $\dD_k$ to $\eqref{4.39}$, we can get
\be\label{4.40}
\begin{aligned}
&\p_t\th_k+u\cdot\na \th_k +\al\Dl u_k+\na\na\cdot u_k=\dD_k\t{F}_1-\al\dD_k\t{F}_3+\mathcal{R}_k,
\end{aligned}
\ee
where $\mathcal{R}_k=[u\cdot\na, \dD_k]\th$. The last two terms on left-hand side of \eqref{4.40} can be written as
\be\label{4.35}
\al\Dl u_k+\na\na\cdot u_k=\al\Dl(\mP u_k+\mvP u_k)+\Dl \mvP u_k=\al\Dl\mP u_k+(1+\al)\Dl\mvP u_k.
\ee
Now inserting \eqref{e4.48} into \eqref{4.35} and substituting the resulting equation into \eqref{4.40}, we can get
\be\label{4.36}
\begin{aligned}
&\p_t\th_k+u\cdot\na \th_k +\f{\al}{\mu_0}\mP \th_k+\f{1+\al}{\la_0+2\mu_0}\mvP \th_k\\
&\qq\qq=-\al\Dl w-(1+\al)\Dl{w}^\bot+\dD_k\t{F}_1-\al\dD_k\t{F}_3+\mathcal{R}_k.
\end{aligned}
\ee
Recalling the fact that $(\mP\th_k|\mvP\th_k)=0$ and $\th_k=\mP\th_k+\mvP\th_k$. A routine procedure shows that after
multiplying $\eqref{4.36}$ by $\mP\th_k$ and $\mvP\th_k$, respectively,
\bea\label{4.59}
&&\|\th_k(t)\|_{L^2}+\int^t_0\|\th_k\|_{L^2}d\tau\ls \|\th_k(0)\|_{L^2}+\int^t_0\|\na u\|_{L^\i}\|\th_k\|_{L^2}d\tau\nn\\
&&\hspace{+10mm}+\int^t_0\|(\La^2 w_k+\La^2{w}^\bot_k)\|_{L^2}d\tau+\int^t_0\|(\dD_k\t{F}_1,\dD_k\t{F}_3,\mathcal{R}_k)\|_{L^2}d\tau.
\eea
Employing commutator estimates in \cite{BCD1} again enable us to get
\be
\sum\limits_{k\in\bZ}2^{ks}\|\mathcal{R}_k\|_{L^2}\ls \|\na u\|_{\dB^{n/2}_{2,1}}\|\th\|_{\dB^{s}_{2,1}}.\nn
\ee
Consequently, multiplying \eqref{4.59} by $2^{k(\f{n}{2}-1)}$ and then summing over the index $k$ satisfying
$2^k>R_0$, we are led to
\bea\label{4.60}
\|\th\|^h_{\t{L}^\i_T\dB^{n/2-1}_{2,1}\cap{\t{L}}^1_T\dB^{n/2-1}_{2,1}}&\ls& \|\th(0)\|^h_{\dB^{n/2-1}_{2,1}}+\|\na u\|_{\t{L}^1_T\dB^{n/2}_{2,1}}\|\th\|_{\t{L}^\i_T\dB^{n/2-1}_{2,1}}\nn\\
&&+\|(w,{w}^\bot)\|_{\t{L}^1_T\dB^{n/2+1}_{2,1}}+\|(\t{F}_1,\t{F}_3)\|^h_{\t{L}^1_T\dB^{n/2-1}_{2,1}}.
\eea
Multiply \eqref{4.60} by a constant $\dl>0$ and then add the resulting inequality to \eqref{4.52} together. By choosing $R_0$ sufficiently large and $\dl>0$ suitably small, we arrive at
 \bea\label{4.600}
&&\|\th\|^h_{\t{L}^\i_T\dB^{n/2-1}_{2,1}\cap{L}^1_T\dB^{n/2-1}_{2,1}}+\|(w,{w}^\bot)\|^h_{\t{L}^\i_T\dB^{n/2-1}_{p,1}\cap L^1_T\dB^{n/2+1}_{p,1}}\nn\\
&\ls&\|(p,\tau)(0)\|^h_{\dB^{n/2}_{2,1}}+\|(w,{w}^\bot)(0)\|^h_{\dB^{n/2-1}_{2,1}}+\|\na u\|_{\t{L}^1_T\dB^{n/2}_{2,1}}\|\th\|_{\t{L}^\i_T\dB^{n/2-1}_{2,1}}\nn\\
&&\q +\|u\cdot\na\th\|^h_{L^1(\dB^{n/2-3}_{2,1})}+\|(\t{F}_1,F_2,\t{F}_3)\|^h_{L^1_T\dB^{n/2-1}_{2,1}}.
 \eea
Clearly, the third term  on the right-hand side of \eqref{4.600} is easily bounded by $\|(p,\tau;u,\th)\|^2_{\E_T}$. The fourth term can be estimated as
\bea
\|u\cdot\na\th\|^h_{L^1(\dB^{n/2-3}_{2,1})}\ls\|u\cdot\na\th\|^h_{L^1(\dB^{n/2-2}_{2,1})}
\ls \|u\|_{\t{L}^2(\dB^{n/2}_{2,1})}\|\na\th\|_{\t{L}^2(\dB^{n/2-2}_{2,1})}
\ls\|(p,\tau;u,\th)\|^2_{\E_T}.
\eea
Finally, keep in mind \eqref{e4.48}, we can conclude that
\bea\label{4.61}
&&\|\th\|^h_{\t{L}^\i_T\dB^{n/2-1}_{2,1}\cap{L}^1_T\dB^{n/2-1}_{2,1}}+\|u\|^h_{\t{L}^\i_T\dB^{n/2-1}_{2,1}\cap L^1_T\dB^{n/2+1}_{2,1}}\nn\\
&\ls&\|(p,\tau)(0)\|^h_{\dB^{n/2}_{2,1}}+\|u(0)\|^h_{\dB^{n/2-1}_{2,1}}+\|(p,\tau;u)\|^2_{\E_T}+\|(\t{F}_1,F_2,\t{F}_3)\|^h_{L^1_T\dB^{n/2-1}_{2,1}}.
 \eea

\subsection{$\boldsymbol{L^\i}$ estimate of $(\boldsymbol{p},\boldsymbol{\tau})$}
In this part, we see that the effective flux $\mathcal{G}$ mentioned plays a key role in deducing the $L^\i$ estimate of $(p,\tau)$ at low frequencies.\\

\noindent\textbf{Step 1: Low-frequency estimates ($\boldsymbol{2^k> R_0}$) }\\
Apply $\dD_k$ to \eqref{r4} to get
\be\label{4.14}
 \left\{
 \begin{aligned}
 &\p_t p_k+u\cdot\na p_k+\nabla\cdot u_k=\dD_k F_1+\mathcal{R}^3_k, \\
 &\p_tu_k+(1-\al)\na p_k-\al\na\cdot\mathcal{G}_k-\mathcal{A} v_k=\dD_kF_2,\\
 &\p_t\mG_k+u\cdot\na \mathcal{G}_k-2D(u_k)=\dD_k(F_3-F_1\text{Id})+\mathcal{R}^5_k,
 \end{aligned}
 \right.
\ee
where $\mathcal{R}^5_k=[u\cdot\na,\dD_k]\mG$.

By taking $L^2$ inner product of $\eqref{4.14}_3$ with $\mG_k$, $\eqref{4.14}_2$ with $\f{2}{\al} u_k$, $\eqref{4.14}_1$ with $\f{2(1-\al)}{\al}p_k$ and then
summing up the resulting equations, we arrive at
\bea\label{4.15}
 &&\f{1}{2}\f{d}{dt}\Big(\|\mG_k\|^2_{L^2}+\f{2}{\al}\|u_k\|^2_{L^2}+\f{2(1-\al)}{\al}\|p_k\|^2_{L^2}\Big)-\f{2}{\al}(\mathcal{A}u_k|u_k)\nn\\
&=&\f{2(1-\al)}{\al}(p_k|\dD_k F_1+\mathcal{R}^3_k)+\f{2}{\al}(u_k|\dD_k F_2)\nn\\
&&\q +\big(\mG_k|\dD_k (F_3-F_1\text{Id})+\mathcal{R}^5_k\big)+\f{1}{2}(|p_k|^2|\na\cdot u)+\f{1}{2}(|\mathcal{G}_k|^2|\na\cdot u)\nn\\
&\triangleq&\hat F(t).
\eea

 Notice that  the coefficient of $\|p_k\|^2_{L^2}$ might be non-positive if $\al\geq1$. In that case, we need to give an auxiliary estimate. Set $\O^\bot\triangleq \La^{-1}\mvP\na\cdot\tau$. Applying $\dD_k$ to $\eqref{r1}_1$, $\dD_k\mvP$ to $\eqref{r1}_2$ and $\dD_k\La^{-1}\mvP\na\cdot$ to $\eqref{r1}_3$, we have
\be\label{e4.3}
\left\{
\begin{aligned}
&\p_t p_k+u\cdot\na p_k+\nabla\cdot \Psi^\bot_k=\dD_k F_1+\mathcal{R}^3_k, \\
&\p_t\Psi^\bot_k-(\la_0+2\mu_0)\Dl\Psi^\bot_k+\na p_k-\al\La\O^\bot_k=\dD_k\mvP F_2,\\
&\p_t\O^\bot_k+u\cdot\na\O^\bot_k+\La  \Psi^\bot_k=\dD_k\La^{-1}\mvP \na\cdot F_3+\mathcal{R}^4_k,\\
\end{aligned}
\right.
\ee
where $\mathcal{R}^3_k=[u\cdot\na,\dD_k]p,\ \mathcal{R}^4_k=[u\cdot\na,\dD_k\La^{-1}\mvP\na\cdot]\tau$.

Taking $L^2$ inner product of $\eqref{e4.3}_1$ with $p_k$, $\eqref{e4.3}_2$ with $\Psi^\bot_k$, $\eqref{e4.3}_3$ with $\al\O^\bot_k$ and then summing up the resulting equations, we arrive at
\bea\label{e4.4}
&&\f{1}{2}\f{d}{dt}\Big(\|p_k\|^2_{L^2}+\|\Psi^\bot_k\|^2_{L^2}+\al\|\O^\bot_k\|^2_{L^2}\Big)+(\la_0+2\mu_0)\|\La\Psi^\bot_k\|^2_{L^2}\nn\\
&&=(p_k|\dD_k F_1+\mathcal{R}^3_k)+\f{1}{2}(|p_k|^2|\na\cdot u)+(\Psi^\bot_k|\dD_k\mvP F_2)\nn\\
&&\q +\f{\al}{2}(|\O^\bot_k|^2|\na\cdot u)+\al(\O^\bot_k|\La^{-1}\mvP \dD_k\na\cdot F_3+\mathcal{R}^4_k)\nn\\
 &&\triangleq\t{F}(t).
\eea

Now, we multiply a small constant $\nu_1>0$ to \eqref{4.15} and add
the resulting equation to \eqref{e4.4}. Choosing $\nu_1>0$ suitably small such that the coefficient of $\|p_k\|^2_{L^2}$ is positive. Consequently,
we have
\bea \label{e4.18}
 &&\f{d}{dt}\Big(\|p_k\|^2_{L^2}+\|u_k\|^2_{L^2}+\|\tau_k\|^2_{L^2}\Big)\ls |\t{F}(t)|+|\hat{F}(t)|.
\eea
Furthermore, bounding the right-hand side of \eqref{e4.18} by Cauchy-Schwarz inequality leads to the following inequality owing to $2^k\leq R_0$,
\bea
&&\q\f{d}{dt}\big(\|p_k\|_{L^2}+\|u_k\|_{L^2}+\|\tau_k\|_{L^2}\big)\nn\\
&&\ls \|\dD_k (F_1,\La^{-1}\t{F}_1,F_3,\La^{-1}\t{F}_3, F_2)\|_{L^2}\nn\\
&&\qq+\|\mathcal{R}^i_k,p_k\na\cdot u,\tau_k\na\cdot u,u\cdot\na\La u_k,\na\cdot u\La u_k\|_{L^2}.\nn
\eea
Noticing that \eqref{e4.10e}-\eqref{e4.13e}, we deduce that
\bea\label{4.19}
&&\|(p,\tau;u)\|^{\ell}_{\t{L}^\i_T\dB^{\frac{n}{2}-1}_{2,1}}\ls \|(a,\tau;u)(0)\|^{\ell}_{\dB^{\frac{n}{2}-1}_{2,1}}\nn\\
&&\qq+\|(F_1,\La^{-1}\t{F}_1, F_2, \La^{-1}\t{F}_3, F_3)\|^\ell_{\t{L}^1_T\dB^{\frac{n}{2}-1}_{2,1}}+\|(p,\tau;u,\th)\|^2_{\mathcal{E}_T}.
\eea

\noindent\textbf{Step 2: High-frequency estimates ($\boldsymbol{2^k> R_0}$) }

Applying $\dD_k$ to the first and third equations of \eqref{r1} gives
\be\label{e4.47}
\left\{
\begin{aligned}
&\p_t p_k+u\cdot\na p_k+\nabla\cdot u_k=\dD_kF_1+\mathcal{R}^6_k, \\
&\p_t\tau_k+u\cdot\na\tau_k+\na\cdot u_k\text{Id}-2D(u_k)=\dD_k F_3+\mathcal{R}^7_k,
\end{aligned}
\right.
\ee
where $\mathcal{R}^6_k=[u\cdot\na,\dD_k]p$ and $\mathcal{R}^7_k=[u\cdot\na,\dD_k]\tau$.

Multiplying the first equation of $\eqref{e4.47}$ by $p_k$ and the second by $\tau_k$, and then integrating over $\bR^n\times [0,t]$, we obtain
\bea\label{e4.48e}
&&\|(p_k,\tau_k)(t)\|_{L^2}\ls \|(p_k,\tau_k)(0)\|_{L^2}+\int^t_0\|\na u_k\|_{L^2}d\tau\nn\\
&&\hspace{10mm} +\int^t_0\|\na u\|_{L^\i}\|(p_k,\tau_k)\|_{L^2}d\tau+\int^t_0\|(\dD_k{F}_1,\dD_k{F}_3,\mathcal{R}^6_k,\mathcal{R}^7_k)\|_{L^2}d\tau.
\eea
It follows from commutator estimates in \cite{BCD1} that
\be
\sum\limits_{k\in\bZ}2^{ks}\|(\mathcal{R}^6_k,\mathcal{R}^7_k)\|_{L^2}\ls \|\na u\|_{\dB^{n/2}_{2,1}}\|(p,\tau)\|_{\dB^{s}_{2,1}}.\nn
\ee
Now multiplying \eqref{e4.48e} by $2^{k\f{n}{2}}$,  and then summing over the index $k$ satisfying
$2^k>R_0$, we are led to
\bea\label{e4.49e}
 \|(p,\tau)\|^h_{\t{L}^\i_T\dB^{n/2}_{2,1}}&\ls& \|(p,\tau)(0)\|^h_{\dB^{n/2}_{2,1}}+\|\na u\|^h_{\t{L}^1_T\dB^{n/2}_{2,1}}\nn\\
&&+\|\na u\|_{\t{L}^1_T\dB^{n/2}_{2,1}}\|(p,\tau)\|_{\t{L}^\i_T\dB^{n/2}_{2,1}}+\|({F}_1,{F}_3)\|^h_{\t{L}^1_T\dB^{n/2}_{2,1}}.
\eea
By using \eqref{eA.9}, we have
\bea \label{e4.50e}
&&\|({F}_1,{F}_3)\|^h_{\t{L}^1_T\dB^{n/2}_{2,1}}\nn\\
&\ls& \|K(a)\|_{\t{L}^\i_T\dB^{n/2}_{2,1}}\|\na u\|_{\t{L}^1_T\dB^{n/2}_{2,1}}
+ \|\tau\|_{\t{L}^\i_T\dB^{n/2}_{2,1}}\|\na u\|_{\t{L}^1_T\dB^{n/2}_{2,1}}\nn\\
&\ls&   (1+\|(p,\tau;u,\th)\|_{\E_T})^{n/2+1}\|(p,\tau;u,\th)\|^2_{\E_T}.
\eea
So, we get
\bea\label{r3e3.39}
 \|(p,\tau)\|^h_{\t{L}^\i_T\dB^{n/2}_{2,1}}&\ls& \|(p,\tau)(0)\|^h_{\dB^{n/2}_{2,1}}+\|\na u\|^h_{\t{L}^1_T\dB^{n/2}_{2,1}}\nn\\
&&+(1+\|(p,\tau;u,\th)\|_{\E_T})^{n/2+1}\|(p,\tau;u,\th)\|^2_{\E_T}.
\eea

\subsection{Estimate for nonlinear terms}
Finally, we devote ourselves to bound those nonlinear terms, which occur in the first two parts. The
following interpolation inequality is frequently used in our analysis
\be\label{4.1er}
\begin{aligned}
&\|f\|_{\t{L}^2_T\dB^{n/2}_{2,1}}\ls \|f\|^{1/2}_{\t{L}^{\i}_T\dB^{n/2-1}_{2,1}}\|f\|^{1/2}_{{L}^{1}_T\dB^{n/2+1}_{2,1}}.
\end{aligned}
\ee

\noindent\textbf{Step 1: Low-frequency estimates ($\boldsymbol{2^k> R_0}$) }\\
Combining \eqref{e4.14} and \eqref{4.19}, we have
\bea\label{4.20}
&&\|(p,\tau)\|^{\ell}_{\t{L}^\i_T\dB^{\frac{n}{2}-1}_{2,1}}+\|(u,\La^{-1}\th)\|^\ell_{\t{L}^\i_T\dB^{\frac{n}{2}-1}_{2,1}\cap\t{L}^1_T\dB^{\frac{n}{2}+1}_{2,1}}\ls \|(a,\tau;u)(0)\|^{\ell}_{\dB^{\frac{n}{2}-1}_{2,1}}\nn\\
&&\qq\qq +\|(F_1,\La^{-1}\t{F}_1, F_2, \La^{-1}\t{F}_3, F_3)\|^\ell_{\t{L}^1_T\dB^{\frac{n}{2}-1}_{2,1}}+\|(p,\tau;u,\th)\|^2_{\mathcal{E}_T}.
\eea
More precisely, we need to deal with the following  nonlinear terms
\be
\ K(a)\na\cdot u,\ \La^{-1}((\na u)^T\na p),\ Q(\tau,\na u),\ \La^{-1}((\na u)^T\cdot\na)\cdot\tau) \q \text{in}\ \t{F}_{1},\ F_1, \t{F_3},\ F_3,\nn
\ee
and
\be
I(a)\mathcal{A}u,\ u\cdot\na u, I(a)\th,\ \f{1}{1+a}\text{div}\big(2\t{\mu}(a)D(u)+\t{\la}(a)\text{div} u\text{Id}\big)  \q \text{in}\ F_{2}.
\ee

Regarding $K(a)\na\cdot u$, by taking $r_1=1,r_2=\i,\ f=\na\cdot u,\ g=K(a)$ in \eqref{eA.9} and using \eqref{2.2}, we have
\bea
&&\sum\limits_{{2^k\leq R_0}}2^{k(n/2-1)}\|\dD_k(K(a)\na\cdot u)\|_{L^1_TL^2}\nn\\
&\ls&\|\na u\|_{\t{L}^{1}_T\dB^{n/2}_{2,1}}\|K(a)\|_{\t{L}^{\i}_T\dB^{n/2-1}_{2,1}}\nn\\
&\ls&(1+\|p\|_{\t{L}^{\i}_T\dB^{n/2}_{2,1}})^{1+[n/2]}\|p\|_{\t{L}^{\i}_T\dB^{n/2-1}_{2,1}}\| u\|_{\t{L}^{1}_T\dB^{n/2+1}_{2,1}}\nn\\
&\ls&(1+\|(p,\tau;u,\th)\|_{\E_T})^{n/2+1}\|(p,\tau;u,\th)\|^2_{\E_T}.
\eea
The terms $Q(\tau,\na u),\ u\cdot\na u $ can be treated along the same line as $K(a)\na u$ by taking $f=\na u$ and $g=\tau, u$ respectively. Also $I(a)\mathcal{A}u$ can be treated by setting
$r_1=\i,r_2=1,\ f=I(a),\ g=\na^2 u$ in \eqref{eA.9} and using \eqref{2.2}. In order to bound the term $\La^{-1}((\na u)^T\na p)$, we apply \eqref{eA.9} by taking $r_1=1,r_2=\i,\ f=\na u,\ g=\na p$. Then we get
\bea
&&\sum\limits_{{2^k\leq R_0}}2^{k(n/2-1)}\|\La^{-1}((\na u)^T\na p)\|_{L^1_TL^2}\nn\\
&=&\sum\limits_{{2^k\leq R_0}}2^{k(n/2-2)}\|((\na u)^T\na p)\|_{L^1_TL^2}\nn\\
&\ls&\|p\|_{\t{L}^{\i}_T\dB^{n/2-1}_{2,1}}\| u\|_{\t{L}^{1}_T\dB^{n/2+1}_{2,1}}
\ls\|(p,\tau;u,\th)\|^2_{\E_T}.
\eea
The term $ \La^{-1}((\na u\cdot\na)\cdot\tau)$ can be treated along the same line as $\La^{-1}((\na u)^T\na p)$.

 For $\ I(a)\th$, we take $r_1=1,\ r_2=\i,\ f=\th,\ g=I(a)$ in \eqref{eA.8} and using \eqref{2.2}. Then we have
\bea
&&\sum\limits_{{2^k\leq R_0}}2^{k(n/2-1)}\|\dD_k(I(a)\th)\|_{L^1_TL^2}\nn\\
&\ls&\|\th\|_{\t{L}^{1}_T\B^{n/2,n/2-1}_{2,1}}\|I(a)\|_{\t{L}^{\i}_T\B^{n/2-1,n/2}_{2,1}}\nn\\
&\ls&(1+\|p\|_{L^\i_TL^\i})^{n/2+1}\|(p,\tau;u,\th)\|^2_{\E_T}\nn\\
&\ls&(1+\|(p,\tau;u,\th)\|_{\E_T})^{n/2+1}\|(p,\tau;u,\th)\|^2_{\E_T}.
\eea
Next we bound nonlinear terms in $F_2$. Denote
\bea
&&I:=\f{1}{1+a}\text{div}\big(2\t{\mu}(a)D(u)\big)\nn\\
&&\q =\f{1}{1+a}\t{\mu}(a)\na^2 u+\f{1}{1+a}\na\t{\mu}(a)\na u\nn\\
&&\q =\t{\mu}(a)\na^2 u-I(a)\t{\mu}(a)\na^2 u+\na\t{\mu}(a)\na u-I(a)\na\t{\mu}(a)\na u\nn\\
&&\q :=I_1+I_2+I_3+I_4.\nn
\eea
The term $I_1$ can treated along the same line as $I(a)\mathcal{A}u$ and $I_3$ can be dealt with by applying \eqref{eA.9} with $f=\na u,\ g=\na \t{\mu}(a)$ and $r_1=1,r_2=\i$.
To bound $I_2$, we have
\bea\label{4.23er}
&&\sum\limits_{{2^k\leq R_0}}2^{k(n/2-1)}\|\dD_k(I(a)\t{\mu}(a)\na^2 v)\|_{L^1_TL^2}\nn\\
&\ls&\|I(a)\|_{\t{L}^{\i}_T\dB^{n/2}_{2,1}}\|\t{\mu}(a)\na^2 u\|_{\t{L}^{1}_T\dB^{n/2-1}_{2,1}}\nn\\
&\ls& \|I(a)\|_{\t{L}^{\i}_T\dB^{n/2}_{2,1}}\|\t{\mu}(a)\|_{\t{L}^{\i}_T\dB^{n/2}_{2,1}}\|\na^2 u\|_{\t{L}^{1}_T\dB^{n/2-1}_{2,1}}\nn\\
&\ls& (1+\|(p,\tau;u,\th)\|_{\E_T})^{n+3}\|(p,\tau;u,\th)\|^2_{\E_T}.
\eea
Regarding $I_4$, it is easy to show that
\bea\label{4.23er}
&&\sum\limits_{{2^k\leq R_0}}2^{k(n/2-1)}\|\dD_k(I(a)\na\t{\mu}(a)\na u)\|_{L^1_TL^2}\nn\\
&\ls&\|I(a)\|_{\t{L}^{\i}_T\dB^{n/2}_{2,1}}\|\na\t{\mu}(a)\na u\|_{\t{L}^{1}_T\dB^{n/2-1}_{2,1}}\nn\\
&\ls& \|I(a)\|_{\t{L}^{\i}_T\dB^{n/2}_{2,1}}\|\na u\|_{\t{L}^{1}_T\dB^{n/2}_{2,1}}\|\na\t{\mu}(a)\|_{\t{L}^{\i}_T\dB^{n/2-1}_{2,1}}\nn\\
&\ls& (1+\|(p,\tau;u,\th)\|_{\E_T})^{n+3}\|(p,\tau;u,\th)\|^2_{\E_T}.
\eea

Since bounding $\f{1}{1+a}\text{div}\big(\t{\la}(a)\text{div} u\text{Id}\big)$ is the same as $I$, we feel free to omit those details. Summing up above all estimates, we conclude that

\bea\label{4.27er}
&&\q\|(p,\tau)\|^{\ell}_{\t{L}^\i_T\dB^{\frac{n}{2}-1}_{2,1}}
+\|(u,\La^{-1}\th)\|^\ell_{\t{L}^\i_T\dB^{\frac{n}{2}-1}_{2,1}\cap\t{L}^1_T\dB^{\frac{n}{2}+1}_{2,1}}\nn\\
&\ls& \|(p,\tau;u)(0)\|^{\ell}_{\dB^{\frac{n}{2}-1}_{2,1}}+(1+\|(p,\tau;u,\th)\|_{\E_T})^{n+3}\|(p,\tau;u,\th)\|^2_{\E_T}.
\eea
\noindent\textbf{Step 2: High-frequency estimates ($\boldsymbol{2^k> R_0}$) }\\
Multiply \eqref{r3e3.39} by a small constant $\nu_2$ and then add the resulting equation to \eqref{4.61}. Note that the term $\|\na u\|^h_{\t{L}^1_T\dB^{n/2}_{2,1}}$ on the right-hand side of \eqref{r3e3.39} can be absorbed by the dissipative term on left-hand side of
 \eqref{4.61}. Consequently, we obtain
\bea
&&\|(p,\tau)\|^h_{\t{L}^\i_T\dB^{n/2}_{2,1}}+\|u\|^h_{\t{L}^\i_T\dB^{n/2-1}_{2,1}\cap L^1_T\dB^{n/2+1}_{2,1}}
+\|\th\|^h_{\t{L}^\i_T\dB^{n/2-1}_{2,1}\cap L^1_T\dB^{n/2-1}_{2,1}}\nn\\
&\ls&\|(p,\tau)(0)\|^h_{\dB^{n/2}_{2,1}}+\|u(0)\|^h_{\dB^{n/2-1}_{2,1}}\nn\\
&&+(1+\|(p,\tau;u,\th)\|_{\E_T})^{n/2+1}\|(p,\tau;u,\th)\|^2_{\E_T}+\|(\t{F}_1,F_2,\t{F}_3)\|^h_{L^1_T\dB^{n/2-1}_{2,1}}.
 \eea
Likely, we bound those nonlinear terms arising in $\t{F}_1,\ F_2,\ \t{F}_3$, see following:
\be
 \na(K(a)\na\cdot u),\ (\na u)^T\na p  \q \text{in}\ \t{F}_{1},\nn
\ee
\be
 \na\cdot Q,\ (\na u\cdot\na)\cdot \tau  \q \text{in}\ \t{F}_3,\nn
\ee
and
\be
I(a)\mathcal{A}u,\ u\cdot\na u, I(a)\th,\ \f{1}{1+a}\text{div}\big(2\t{\mu}(a)D(u)+\t{\la}(a)\text{div} u\text{Id}\big)  \q \text{in}\ F_{2}.
\ee
In order to bound $(\na u)^T\na p$, by \eqref{eA.9}, we have \bea\label{4.48rt}
&&\sum\limits_{{2^k> R_0}}2^{k(n/2-1)}\|\dD_k\big((\na u)^T\na p\big)\|_{L^1_TL^2}\nn\\
&\ls&\|\na u\|_{\t{L}^{1}_T\dB^{n/2}_{2,1}}\|\na p\|_{\t{L}^{\i}_T\dB^{n/2-1}_{2,1}}
\ls\| u\|_{\t{L}^{1}_T\dB^{n/2+1}_{2,1}}\| p\|_{\t{L}^{\i}_T\dB^{n/2}_{2,1}}  \nn\\
&\ls&\|(p,\tau;u,\th)\|^2_{\E_T}.
\eea

Regarding $\na(K(a)\na\cdot u)$, we write $\na(K(a)\na\cdot u)=K(a)\na\na\cdot u+\na\cdot u \na K(a)$. The estimate for $\na\cdot u \na K(a)$ can be handled with at the same way as $(\na u)^T\na p$. For $K(a)\na\na\cdot u$, we arrive at
\bea
&&\sum\limits_{{2^k> R_0}}2^{k(n/2-1)}\|\dD_k\big(K(a)\na\na\cdot u\big)\|_{L^1_TL^2}\nn\\
&\ls&\|K(a)\|_{\t{L}^{\i}_T\dB^{n/2}_{2,1}}\|\na^2 u\|_{\t{L}^{1}_T\dB^{n/2-1}_{2,1}}\nn\\
&\ls&(1+\|(p,\tau;u,\th)\|_{\E_T})^{n/2+1}\|(p,\tau;u,\th)\|^2_{\E_T}.
\eea

Bounding $\t{F}_3$ can be treated along the same line as $\t{F}_1$. The high frequency of $F_2$ can be dealt with at the similar way as the low frequency, which is left to the interested reader. Consequently, we deduce that
\bea\label{e4.47e}
 &&\|\th\|^h_{\t{L}^\i_T\dB^{n/2-1}_{2,1}\cap{L}^1_T\dB^{n/2-1}_{2,1}}+\|u\|^h_{\t{L}^\i_T\dB^{n/2-1}_{2,1}\cap L^1_T\dB^{n/2+1}_{2,1}}\nn\\
&\ls&\|(p,\tau)(0)\|^h_{\dB^{n/2}_{2,1}}+\|u(0)\|^h_{\dB^{n/2-1}_{2,1}}+(1+\|(p,\tau;u,\th)\|_{\E_T})^{n+3}\|(p,\tau;u)\|^2_{\E_T}.
 \eea
At last, combining \eqref{r3e3.39} and \eqref{e4.47e}, we achieve the high-frequency estimate
 \bea\label{4.51er}
&&\|(p,\tau)\|^h_{\t{L}^\i_T\dB^{n/2}_{2,1}}+\|u\|^h_{\t{L}^\i_T\dB^{n/2-1}_{2,1}\cap L^1_T\dB^{n/2+1}_{2,1}}
+\|\th\|^h_{\t{L}^\i_T\dB^{n/2-1}_{2,1}\cap L^1_T\dB^{n/2-1}_{2,1}}\nn\\
&\ls&\|(p_0,\tau_0)\|^h_{\dB^{n/2}_{2,1}}+\|u_0\|^h_{\dB^{n/2-1}_{2,1}}+(1+\|(p,\tau;u,\th)\|_{\E_T})^{n+3}\|(p,\tau;u,\th)\|^2_{\E_T}.
 \eea

\noindent The inequality (\ref{4.1}) is followed by \eqref{4.27er} and \eqref{4.51er}. Hence, the proof of Proposition \ref{proapriori} is  complete. \ef

\section{Proof of Theorem \ref{thglobal1}}
Let us recall a local-in-time existence result of \eqref{R-E1}-\eqref{id0} which has been achieved in \cite{QZ1}.

\begin{proposition}\label{prop4.1}
Assume $(\rho_0-1,F_0-I)\in \big(\dot{B}^{n/2}_{2,1}\big)^{1+n^2}$ and $u_0\in\big(\dot{B}^{n/2-1}_{2,1}\big)^n$ with $\rho_0$ bounded away from $0$. There exists a time $T>0$ such that \eqref{R-E1}-\eqref{id0} has a unique solution $(\rho,F;u)$ with $\rho$ bounded away from zero and
\be
(\rho-1,F-I)\in \Big(C([0,T);\dot{B}^{n/2}_{2,1})\Big)^{1+n^2}, u\in \Big(C([0,T);\dot{B}^{n/2-1}_{2,1})\cap L^1([0,T);\dot{B}^{n/2+1}_{2,1})\Big)^n
.\nn
\ee \end{proposition}

Based on Proposition \ref{prop4.1}, the proof of Theorem \ref{thglobal1} can be finished by the standard continuity argument. Indeed, Proposition \ref{prop4.1} indicates that there exists a maximal time $T>0$ such that system \eqref{R-E1} admits a unique solution. Clearly, the system \eqref{r2} also has a solution $(p,\tau;u)$ which locally exits on $[0,T)$. It follows from the assumption of Theorem \ref{thglobal1} and Lemma \ref{l2.3} that
\be
\|(p_{0},\tau_{0};u_{0})\|_{\mathcal{E}_0}\leq C_0\eta,\nn
\ee
for some positive constant $C_0$. Fixed a constant $M>0$ (to be determined later), we define
\be
T^\ast \triangleq \sup\{t\in[0,T)\big| \|(p,\tau;u,\th)\|_{\mathcal{E}_t}\leq M\eta\}.\nn
\ee
 Claim that
\be
T^\ast=T.\nn
\ee
According to the continuity argument, it suffices to show \be\label{4.15x}
\|(p,\tau;u,\th)\|_{\mathcal{E}_{T}}\leq \f{1}{2}M\eta.
\ee
Indeed, noting that
\be
\|a\|_{L^\i([0,T)\times\bR^n)}=\|h(p)\|_{L^\i([0,T)\times\bR^n)}\leq C_1\|p\|_{L^\i_{T}{\B^{n/2-1,n/2}}}.\nn
\ee
We can choose $\eta$ sufficiently small such that
\be
M\eta\leq \f{1}{\eta_0C_1},\nn
\ee
so
\be
\|a\|_{L^\i([0,T)\times\bR^n)}\leq \eta_0.\nn
\ee
By applying Proposition \ref{proapriori}, we obtain
\begin{equation}\label{4.155}
\|(p,\tau;u,\th)\|_{\mathcal{E}_{T}}\leq C\{C_0\eta+(M\eta)^{2}(1+M\eta)^{n+3}\}.
\end{equation}
By choosing $M=3CC_0$ and $\eta$ sufficient small enough such that
\be
 C(M\eta)(1+M\eta)^{n+3}\leq \f{1}{6},\nn
\ee
so \eqref{4.15x} is followed by \eqref{4.155} directly.
Actually the above argument implies
\be\label{4.55}
\|(p,\tau;u,\th)\|_{\E_{T}}\leq C\|(p,\tau;u)\|_{\E_0}. \ee
 Consequently, the continuity argument ensures that $T=+\i$. It follows from
the third equation of \eqref{R-E2} that
 \be
\p_t(F-I)+u\cdot\nabla (F-I)=\na u+\nabla u(F-I). \nn
 \ee
By using Lemma \ref{ll2.5}, we have
\bea
&&\|F-I\|_{\t{L}_T^{\infty}\dB^{n/2}_{2,1}}\nn\\
&\leq& \exp(C\int^\i_0\|\na u\|_{\dB^{n/2}_{2,1}}d\tau)\big(\|F_0-I\|_{\dB^{n/2}_{2,1}}
+\int^\i_0(\|\nabla u\|_{\dB^{n/2}_{2,1}}+\|\nabla u(F-I)\|_{\dB^{n/2}_{2,1}})d\tau\big)\nn\\
&\leq& C\big(\|F_0-I\|_{\dB^{n/2}_{2,1}}+\|F-I\|_{\t{L}^{\infty}_T\dB^{n/2}_{2,1}}\|\na u\|_{\t{L}^1_T\dB^{n/2}_{2,1}}+M\eta\big)\nn\\
&\leq&   C\|F_0-I\|_{\dB^{n/2}_{2,1}}+CM\eta\|F-I\|_{\t{L}_T\dB^{n/2}_{2,1}}+CM\eta. \nn
\eea
Furthermore, we chose $\eta$ small enough such that $CM\eta\leq 1/2$ and thus obtain
\be\label{4.56}
 \|F-I\|_{\t{L}_T^{\infty}\dB^{n/2}_{2,1}}\leq C \|F_0-I\|_{\dB^{n/2}_{2,1}}+CM\eta.
\ee
The continuity argument and \eqref{4.55}-\eqref{4.56} enable us to finish the proof of Theorem \ref{thglobal1} eventually.

\vskip 1cm

\section{Appendix: analysis tools}
To make the manuscript self-contained as soon as possible, we would like to collect nonlinear estimates in the last section. See \cite{BCD1} for more details.

\begin{lemma}\label{l2.2}

For the Besov space, we have the following properties:
\begin{itemize}
 \item $\B^{s_2,\sigma}\subseteq \B^{s_1,\sigma}$ for $s_1\geq s_2$ and $\B^{s,\sigma_2}\subseteq \B^{s,\sigma_1}$ for $\sigma_1\leq \sigma_2$.

 \item Interpolation: For $s_1,s_2,\sigma_1,\sigma_2\in \bR$ and $\theta\in [0,1]$, we have
\be
\|f\|_{\B^{\theta s_1+(1-\theta)s_2,\theta\sigma_1+(1-\theta)\sigma_2}}\leq \|f\|^\theta_{\B^{s_1,\sigma_1}}\|f\|^{(1-\theta)}_{\B^{s_2,\sigma_2}}.\nn
\ee
\end{itemize}
\end{lemma}
System (1.3) involves compositions of functions and they are bounded according
to the following lemma.
\begin{lemma} \label{l2.3}
Let $F:\R\rightarrow\R$ be  smooth
with $F(0)=0.$
For  all  $1\leq p,r\leq\infty$ and  $s>0$, we have
\be\label{2.2}
\|F(f)\|_{\tilde{L}^{r}_{T}(\dot B^s_{p,1})}\leq C\|f\|_{\tilde{L}^{r}_{T}(\dot B^s_{p,1})},
\ee
where  $C$ depending only on $\|f\|_{L_{T}(L^\infty)},$ $F'$ (and higher derivatives),  $s,$ $p$ and $n.$
\end{lemma}

For the heat equation, one has the following optimal regularity estimate.
\begin{lemma}\label{lheat}
  Let $p,r\in [1,\i]$, $s\in\bR$, and $1\leq \rho_2\leq \rho_1\leq\i $ Assume that $u_0\in\dB^{s-1}_{p,r}$, $f\in \t{L}^{\rho_2}_T\dB^{s-3+\f{2}{\rho_2}}_{p,r}$. Let $u$ be a solution of the equation
 \be
 \p_tu -\mu \Dl u=f,\q u|_{t=0}=u_0.\nn
 \ee
 Then for $t\in[0,T]$, there holds
 \be\label{2.3}
 \mu^{\f{1}{\rho_1}}\|u\|_{\t{L}^{\rho_1}_T\dB^{s-1+2/\rho_1}_{p,r}}\leq C\big(\|u_0\|_{\dB^{s-1}_{p,r}}+\mu^{1/\rho_2-1}\|f\|_{\t{L}^{\rho_2}_T\dB^{s-3+\f{2}{\rho_2}}_{p,r}}\big).
 \ee
\end{lemma}

In order to obtain the $L^{\infty}$ estimate of original variable $F$ with respect to time $t$, we need the estimate for the transport equation.

 \begin{lemma}\label{ll2.5}
  Let $s\in(-n\min(1/p,1/p'),1+n/p)$ and $1\leq p,q\leq \i$. Let $v$ be a vector field such that $\na v\in L^1_T\dB^{n/p}_{p,1}$. Assume that $f_0\in\dB^s_{p,q},g\in L^1_T\dB^s_{p,q}$, and $f$ is a solution of the transport equation
 \be
 \p_tf+v\cdot\na f=g,\q f|_{t=0}=f_0. \nn
 \ee
 Then for $t\in[0,T]$, it holds that
 \be
 \|f\|_{\t{L}_t\dB^s_{p,q}}\leq \exp\big(C\int^t_0\|\na v(\tau)\|_{\dB^{n/p}_{p,1}}d\tau\big)\big(\|f_0\|_{\dB^s_{p,q}}+\int^t_0\|g(\tau)\|_{\dB^s_{p,q}}d\tau\big).\nn
 \ee
 \end{lemma}

The standard product estimate is also used in our analysis.
\begin{proposition}[\cite{Cjy:2004}]
 Let $1\leq r,r_1,r_2\leq \i$ with $\f{1}{r}=\f{1}{r_1}+\f{1}{r_2}$ and $s,t\leq n/2$, $s+t\geq 0$. Then we have
 \bea\label{eA.9}
  &&\|fg\|_{\t{L}^r_t\dB^{s+t-n/2}_{2,1}}\ls \|f\|_{\t{L}^{r_1}_t\dB^{s}_{2,1}}\|g\|_{\t{L}^{r_2}_t\dB^{t}_{2,1}}.
  \eea
\end{proposition}

In addition, we develop a product estimate in the framework of hybrid Besov spaces by using Bony's decompositions. Let us by denote
 $\chi\{\cdot\}$ the characteristic function in $\bZ$ and $\{c(j)\}_{j\in\bZ}$ be some sequence on $\ell^1$ satisfying $\|\{c(j)\}\|_{\ell^1}=1$.
\begin{lemma}\label{lA.2}
 Let $s,t,\sigma, \tau\in\bR$. Then we have the following:

(i) For $2^j\leq R_0$, if $s \leq n/2$, then
\bea\label{eA.1}
\|\dD_j(T_fg)\|_{L^2}\leq Cc(j)2^{j(n/2-s-t)}\|f\|^\ell_{\dB^s_{2,1}}\|g\|^\ell_{\dB^t_{2,1}}.
\eea

(ii) For $2^j> R_0$, if $s,\sigma \leq n/2$, then
\bea\label{eA.2}
\|\dD_j(T_fg)\|_{L^2}
\leq Cc(j)\big(2^{j(n/2-s-t)}\|f\|^\ell_{\dB^s_{2,1}}\|g\|^h_{\dB^t_{2,1}}
+2^{j(n/2-\sigma-\tau)}\|f\|^h_{\dB^\sigma_{2,1}}\|g\|^h_{\dB^\tau_{2,1}}\big).
\eea
\end{lemma}
\pf  With our choice of $\varphi$ in Introduction, it is easy to see that
\be\label{ee}
\begin{aligned}
&\dD_j\dD_kf=0\q \text{if}\ |j-k|\geq 2,\\
&\dD_j(\dot{S}_{k-1}f\dD_kf)=0\q \text{if}\ |j-k|\geq 5.
\end{aligned}
\ee
Thanks to \eqref{ee}, we have
\bea
\dD_j(T_fg)&=&\sum\limits_{|k-j|\leq 4}\dD_j(\dot{S}_{k-1}f\dD_kg)\nn\\
&=&\sum\limits_{|k-j|\leq 4}\sum_{k'\leq k-2}\dD_j(\dD_{k'}f\dD_kg).\nn
\eea
Denote $J:=\{(k,k'):|k-j|\leq 4,k'\leq k-2\}$, then for $2^j\leq R_0$,
\bea
\|\dD_j(T_fg)\|_{L^2}&\leq& \sum\limits_{J}\|\dD_j(\dD_{k'}f\dD_kg)\|_{L^2}\nn\\
&\ls&\sum\limits_{J}2^{k's}\|\dD_{k'}f\|_{L^2}2^{k'(n/2-s)}2^{kt}\|\dD_{k}g\|_{L^2}2^{-kt}\nn\\
&\ls&c(j)2^{j(n/2-s-t)}\|f\|^\ell_{\dB^{s}_{2,1}}\|g\|^\ell_{\dB^{t}_{2,1}}.\nn
\eea
Next we turn to prove \eqref{eA.2}. We write $J=J_1+J_2$, where
 \be
 J_1=J\cap \{2^{k'}\leq R_0\},\q J_2=J\cap \{2^{k'}> R_0\}.\nn
 \ee
  For $2^j> R_0$ and $s,\sigma\leq n/2$, one has
\bea
&&\|\dD_j(T_fg)\|_{L^2}\nn\\
&\leq& \sum\limits_{J}\|\dD_j(\dD_{k'}f\dD_kg)\|_{L^2}\nn\\
&\ls&\sum\limits_{J}\|\dD_{k'}f\|_{L^\i}\|\dD_{k}g\|_{L^2}\nn\\
&\ls&\sum\limits_{J_1}\|\dD_{k'}f\|_{L^\i}\|\dD_{k}g\|_{L^2}+\sum\limits_{J_2}\|\dD_{k'}f\|_{L^\i}\|\dD_{k}g\|_{L^2}\nn\\\nn\\
&\ls&\sum\limits_{J_1}2^{k's}\|\dD_{k'}f\|_{L^2}2^{k'(n/2-s)}2^{kt}\|\dD_{k}g\|_{L^2}2^{-kt}\nn\\
&& +\sum\limits_{J_1}2^{k'\sigma}\|\dD_{k'}f\|_{L^2}2^{k'(n/2-\sigma)}2^{k\tau}\|\dD_{k}g\|_{L^2}2^{-k\tau}\nn\\
&\ls&c(j)2^{j(n/2-s-t)}\|f\|^\ell_{\dB^{s}_{2,1}}\|g\|^h_{\dB^{t}_{2,1}}+c(j)2^{j(n/2-\sigma-\tau)}\|f\|^h_{\dB^{\sigma}_{2,1}}\|g\|^h_{\dB^{\tau}_{2,1}},\nn
\eea which is just \eqref{eA.2}.
\ef

\begin{lemma}\label{lA.2}
 Let $s,t,\sigma,\tau\in\bR$. Assume that $s+t\geq0,\sigma+\tau\geq 0$. It holds that
 \bea\label{eA.3}
  &&\|\dD_jR(f,g)\|_{L^2}\nn\\
&\leq&Cc(j)\big(2^{j(n/2-s-t)}\|f\|^\ell_{\dB^s_{2,1}}\|g\|^\ell_{\dB^t_{2,1}}+2^{j(n/2-\sigma-\tau)}\|f\|^h_{\dB^\sigma_{2,1}}\|g\|^h_{\dB^\tau_{2,1}}\big).
 \eea
\end{lemma}
\pf Thanks to \eqref{ee}, we have
\be
\dD_jR(f,g)=\sum\limits_{k\geq j-3}\sum\limits_{|k-k'|\leq 1}\dD_j(\dD_kf\dD_{k'}g).\nn
\ee

Denote $J:=\{(k,k'):k\geq j-3,|k-k'|\leq 1\}$ and
\be
J_1=J\cap \{2^{k'}\leq R_0\},\q J_2=J\cap \{2^{k'}> R_0\}.\nn
\ee
 Then when $s+\tau\geq 0$, we have
\bea
&&\|\dD_jR(f,g)\|_{L^2}\nn\\
&\ls& 2^{jn/2}\sum\limits_{(k,k')\in J}\|\dD_kf\dD_{k'}g\|_{L^{1}}\nn\\
&=& 2^{jn/2}\sum\limits_{(k,k')\in J_1}\|\dD_kf\dD_{k'}g\|_{L^{1}}+2^{jn/2}\sum\limits_{(k,k')\in J_2}\|\dD_kf\dD_{k'}g\|_{L^{1}}\nn\\
&\ls&  2^{jn/2}\sum\limits_{(k,k')\in J_1}\|\dD_kf\|_{L^{2}}\|\dD_{k'}g\|_{L^{2}}+2^{jn/2}\sum\limits_{(k,k')\in J_2}\|\dD_kf\|_{L^{2}}\|\dD_{k'}g\|_{L^{2}}\nn\\
&\ls&  2^{jn/2}\sum\limits_{(k,k')\in J_1}2^{ks}\|\dD_kf\|_{L^{2}}2^{-ks}2^{k't}\|\dD_{k'}g\|_{L^{2}}2^{-k't}\nn\\
&&+ 2^{jn/2}\sum\limits_{(k,k')\in J_2}2^{k\sigma}\|\dD_kf\|_{L^{2}}2^{-k\sigma}2^{k'\tau}\|\dD_{k'}g\|_{L^{2}}2^{-k'\tau}\nn\\
&\ls&  c(j)2^{j(n/2-s-t)}\|f\|^\ell_{\dB^s_{2,1}}\|g\|^\ell_{\t{L}^{r_2}_t\dB^t_{2,1}}+c(j)2^{j(n/2-\sigma-\tau)}\|f\|_{\dB^\sigma_{2,1}}\|g\|_{\t{L}^{r_2}_t\dB^\tau_{2,1}}.\nn
\eea

This finishes the proof of Lemma \ref{lA.2}.\ef

Having above continuity of para-product operator and remaining operator, one can get the key product estimate.
\begin{proposition}\label{proparaproduct}
It holds that
  \bea\label{eA.4}
  &&\|fg\|_{\dB^{n/2-1}_{2,1}}\ls \|f\|_{\B^{n/2,n/2-1}_{2,1}}\|g\|_{\B^{n/2-1,n/2}_{2,1}}.
  \eea
\end{proposition}
\pf By Bony decomposition, we write $fg=T_fg+T_gf+R(f,g)$. At low frequencies, we take $s=n/2,t=n/2-1$  in
\eqref{eA.1} for $T_fg$ and $s=n/2-1,t=n/2$ in \eqref{eA.1} for $T_gf$. Then we get
 \be\label{eA.5}
\sum\limits_{2^j\leq R_0}2^{j(n/2-1)}\big(\|\dD_j(T_fg)\|_{L^2}+\|\dD_j(T_gf)\|_{L^2}\big)\ls\|f\|^\ell_{\dB^{n/2}_{2,1}}\|g\|^\ell_{\dB^{n/2-1}_{2,1}}.
 \ee
For the high frequency, we choose $s=n/2,t=n/2-1,\sigma=n/2-1,\tau=n/2$  in
\eqref{eA.2} for $T_fg$; $s=n/2-1,t=n/2,\sigma=n/2,\tau=n/2-1$ in \eqref{eA.2} for $T_gf$, which lead to
 \be\label{eA.6}
\sum\limits_{2^j\geq R_0}2^{j(n/2-1)}\big(\|\dD_j(T_fg)\|_{L^2}+\|\dD_j(T_gf)\|_{L^2}\big)\ls\|f\|_{\B^{n/2,n/2-1}}\|g\|_{\B^{n/2-1,n/2}}.
 \ee
Finally, we choose  $s=n/2,t=n/2-1,\sigma=n/2-1,\tau=n/2$ in \eqref{eA.3} and get
\be\label{eA.7}
\sum\limits_{j\in\bZ}2^{j(n/2-1)}\|\dD_j(R(f,g)\|_{L^2}\ls\|f\|^\ell_{\dB^{n/2}_{2,1}}\|g\|^\ell_{\B^{n/2-1}_{2,1}}
+\|f\|^h_{\dB^{n/2-1}_{2,1}}\|g\|^h_{\dB^{n/2}_{2,1}}.
 \ee
 Combining \eqref{eA.5}, \eqref{eA.6} and \eqref{eA.7}, we arrive at \eqref{eA.4}. Therefore, the proof of Proposition \ref{proparaproduct} is completed.

Finally, let us point out the new product estimate remains valid in Chemin-Lerner's spaces whereas the time exponent $r$ behaves according to the H\"{o}lder inequality.
 \begin{remark} The inequality
\bea\label{eA.8}
  &&\|fg\|_{\t{L}^r_t\dB^{n/2-1}_{2,1}}\ls \|f\|_{\t{L}^{r_1}_t\B^{n/2,n/2-1}_{2,1}}\|g\|_{\t{L}^{r_2}_t\B^{n/2-1,n/2}_{2,1}}
  \eea
 holds whenever $1\leq r,r_1,r_2\leq \i$ and $\f{1}{r}=\f{1}{r_1}+\f{1}{r_2}$.
\end{remark}

\noindent
{\bf Acknowledgement.}
This first author (X. Pan) is supported by Natural Science Foundation of Jiangsu Province (SBK2018041027) and National Natural Science Foundation of China (11801268). The second author (J. Xu) is supported by the National Natural Science Foundation of China (11871274) and supported by the China Scholarship Council (201906835023). The paper is partially finished during his visit at Waseda University. He would like to thank Professors Shuichi Kawashima and Yoshihiro Shibata for their warm hospitality. The third author (Y. Zhu) is supported by the National Natural Science Foundation of China (11801175).

\qq \qq\qq\qq\qq\qq Xinghong Pan

\qq\qq\qq\qq\qq\qq Department of Mathematics,

\qq\qq\qq\qq\qq\qq Nanjing University of Aeronautics and Astronautics,

\qq\qq\qq\qq\qq\qq Nanjing 211106, People's Republic of China.

\qq\qq\qq\qq\qq\qq Email:xinghong\underline{\ \ }87@nuaa.edu.cn

\vskip 0.3cm

\qq \qq\qq\qq\qq\qq Jiang Xu

\qq\qq\qq\qq\qq\qq Department of Mathematics,

\qq\qq\qq\qq\qq\qq Nanjing University of Aeronautics and Astronautics,

\qq\qq\qq\qq\qq\qq Nanjing 211106, People's Republic of China.

\qq\qq\qq\qq\qq\qq Email:jiangxu\underline{\ \ }79@nuaa.edu.cn

\vskip 0.3cm

\qq \qq\qq\qq\qq\qq Yi Zhu

\qq\qq\qq\qq\qq\qq Department of Mathematics,

\qq\qq\qq\qq\qq\qq East China University of Science and Technology,

\qq\qq\qq\qq\qq\qq Shanghai 200237, People's Republic of China.

\qq\qq\qq\qq\qq\qq Email: zhuyim@ecust.edu.cn

 \end{document}